\documentclass[11pt]{amsart}
\usepackage{amsmath}
\usepackage{amsxtra}
\usepackage{amscd}
\usepackage{amsthm}
\usepackage{amsfonts}
\usepackage{amssymb}
\usepackage{eucal}

\newtheorem{cor}{Corollary}
\newtheorem{lem}{Lemma}
\newtheorem{prop}{Proposition}

\newtheorem{thm}{Theorem}

\theoremstyle{definition}

\theoremstyle{remark}

\newcommand{\nc}{\newcommand}
\nc{\renc}{\renewcommand}
\nc{\ssec}{\subsection}
\nc{\sssec}{\subsubsection}
\nc{\on}{\operatorname}

\nc\ol{\overline}
\nc\wt{\widetilde}
\nc\tboxtimes{\wt{\boxtimes}}
\nc{\alp}{\alpha}

\nc{\ZZ}{{\mathbb Z}}
\nc{\NN}{{\mathbb N}}
\nc{\CC}{{\mathbb C}}
\nc{\OO}{{\mathbb O}}
\renc{\SS}{{\mathbb S}}
\nc{\DD}{{\mathbb D}}
\nc{\GG}{{\mathbb G}}
\renewcommand{\AA}{{\mathbb A}}
\nc{\Fq}{{\mathbb F}_q}
\nc{\Fqb}{\ol{{\mathbb F}_q}}
\nc{\Ql}{\ol{{\mathbb Q}_\ell}}
\nc{\id}{\text{id}}
\nc\X{\mathcal X}

\nc{\Hom}{\on{Hom}}
\nc{\Lie}{\on{Lie}}
\nc{\Loc}{\on{Loc}}
\nc{\Pic}{\on{Pic}}
\nc{\Bun}{\on{Bun}}
\nc{\IC}{\on{IC}}
\nc{\Aut}{\on{Aut}}
\nc{\rk}{\on{rk}}
\nc{\Sh}{\on{Sh}}
\nc{\Perv}{\on{Perv}}
\nc{\pos}{{\on{pos}}}
\nc{\Conv}{\on{Conv}}
\nc{\Sph}{\on{Sph}}
\nc{\Sym}{\on{Sym}}
\nc{\BunBb}{\overline{\Bun}_B}
\nc{\Buno}{\overset{o}{\Bun}}
\nc{\BunPb}{{\overline{\Bun}_P}}
\nc{\BunBM}{\overline{\Bun}_{B(M)}}
\nc{\BunPbw}{{\widetilde{\Bun}_P}}
\nc{\BunBP}{\widetilde{\Bun}_{B,P}}
\nc{\GUb}{\overline{G/U}}
\nc{\GUPb}{\overline{G/U(P)}}

\nc{\Hhom}{\underline{\on{Hom}}}
\nc\syminfty{\on{Sym}^{\infty}}
\nc\lal{\ol{\lambda}}
\nc\xl{\ol{x}}
\nc\thl{\ol{\theta}}
\nc\nul{\ol{\nu}}
\nc\mul{\ol{\mu}}
\nc\Sum\Sigma
\nc{\oX}{\overset{o}{X}{}}

\nc{\M}{{\mathcal M}}
\nc{\N}{{\mathcal N}}
\nc{\F}{{\mathcal F}}
\nc{\D}{{\mathcal D}}
\nc{\Q}{{\mathcal Q}}
\nc{\Y}{{\mathcal Y}}
\nc{\G}{{\mathcal G}}
\nc{\E}{{\mathcal E}}
\nc{\CalC}{{\mathcal C}}
\nc\Dh{\widehat{\D}}
\renewcommand{\O}{{\mathcal O}}
\nc{\C}{{\mathcal C}}
\nc{\K}{{\mathcal K}}
\renewcommand{\H}{{\mathcal H}}

\nc{\T}{{\mathcal T}}
\nc{\V}{{\mathcal V}}
\renc{\P}{{\mathcal P}}
\nc{\A}{{\mathcal A}}
\nc{\B}{{\mathcal B}}
\nc{\U}{{\mathcal U}}

\nc{\Gr}{\on{Gr}}

\nc{\frn}{{\check{\mathfrak u}(P)}}
\nc\f{{\mathfrak f}}

\nc{\q}{{\mathfrak q}}
\nc{\p}{{\mathfrak p}}
\nc{\s}{{\mathfrak s}}
\nc\w{\text{w}}

\nc\mathi\iota
\nc\Spec{\on{Spec}}
\nc\Mod{\on{Mod}}
\nc{\tw}{\widetilde{\mathfrak t}}
\nc{\pw}{\widetilde{\mathfrak p}}
\nc{\qw}{\widetilde{\mathfrak q}}
\nc{\jw}{\widetilde j}

\nc{\grb}{\overline{\Gr}}
\nc{\I}{\mathcal I}

\nc{\lambdach}{{\check\lambda}}
\nc{\Lambdach}{{\check\Lambda}{}}
\nc{\much}{{\check\mu}}
\nc{\omegach}{{\check\omega}}
\nc{\nuch}{{\check\nu}}
\nc{\etach}{{\check\eta}}
\nc{\alphach}{{\check\alpha}}
\nc{\betach}{{\check\beta}}
\nc{\rhoch}{{\check\rho}}
\nc{\ch}{{\check h}}

\nc{\Hb}{\overline{\H}}

\nc{\HC}{{\mathcal{HC}}}

\nc{\BA}{{\mathbb{A}}}
\nc{\BC}{{\mathbb{C}}}
\nc{\BM}{{\mathbb{M}}}
\nc{\BN}{{\mathbb{N}}}
\nc{\BP}{{\mathbb{P}}}
\nc{\BR}{{\mathbb{R}}}
\nc{\BZ}{{\mathbb{Z}}}
\nc{\BS}{{\mathbb{S}}}

\nc{\CA}{{\mathcal{A}}}
\nc{\CB}{{\mathcal{B}}}
\nc{\CalD}{{\mathcal D}}
\nc{\CE}{{\mathcal{E}}}
\nc{\CF}{{\mathcal{F}}}
\nc{\CG}{{\mathcal{G}}}
\nc{\CL}{{\mathcal{L}}}
\nc{\CM}{{\mathcal{M}}}
\nc{\CMM}{{\mathcal{M}^{\operatorname{gen}}_\hbar(-\rho)}}
\nc{\CN}{{\mathcal{N}}}
\nc{\CO}{{\mathcal{O}}}
\nc{\CP}{{\mathcal{P}}}
\nc{\CQ}{{\mathcal{Q}}}
\nc{\CR}{{\mathcal{R}}}
\nc{\CS}{{\mathcal{S}}}
\nc{\CT}{{\mathcal{T}}}
\nc{\CU}{{\mathcal{U}}}
\nc{\CV}{{\mathcal{V}}}
\nc{\CW}{{\mathcal{W}}}
\nc{\CZ}{{\mathcal{Z}}}

\nc{\gen}{{\operatorname{gen}}}
\nc{\cM}{{\check{\mathcal M}}{}}
\nc{\csM}{{\check{\mathcal A}}{}}
\nc{\oM}{{\overset{\circ}{\mathcal M}}{}}
\nc{\obM}{{\overset{\circ}{\mathbf M}}{}}
\nc{\oCA}{{\overset{\circ}{\mathcal A}}{}}
\nc{\obA}{{\overset{\circ}{\mathbf A}}{}}
\nc{\ooM}{{\overset{\circ}{M}}{}}
\nc{\osM}{{\overset{\circ}{\mathsf M}}{}}
\nc{\vM}{{\overset{\bullet}{\mathcal M}}{}}
\nc{\nM}{{\underset{\bullet}{\mathcal M}}{}}
\nc{\oD}{{\overset{\circ}{\mathcal D}}{}}
\nc{\obD}{{\overset{\circ}{\mathbf D}}{}}
\nc{\oA}{{\overset{\circ}{\mathbb A}}{}}
\nc{\op}{{\overset{\bullet}{\mathbf p}}{}}
\nc{\cp}{{\overset{\circ}{\mathbf p}}{}}
\nc{\oU}{{\overset{\bullet}{\mathcal U}}{}}
\nc{\oZ}{{\overset{\circ}{\mathcal Z}}{}}
\nc{\ofZ}{{\overset{\circ}{\mathfrak Z}}{}}

\nc{\fa}{{\mathfrak{a}}}
\nc{\fb}{{\mathfrak{b}}}
\nc{\fg}{{\mathfrak{g}}}
\nc{\fgl}{{\mathfrak{gl}}}
\nc{\fh}{{\mathfrak{h}}}
\nc{\fj}{{\mathfrak{j}}}
\nc{\fm}{{\mathfrak{m}}}
\nc{\fn}{{\mathfrak{n}}}
\nc{\fu}{{\mathfrak{u}}}
\nc{\fp}{{\mathfrak{p}}}
\nc{\frr}{{\mathfrak{r}}}
\nc{\fs}{{\mathfrak{s}}}
\nc{\fT}{{\mathfrak{T}}}
\nc{\ofT}{{\overline{\mathfrak T}}}
\nc{\ofS}{{\overline{\mathfrak S}}}
\nc{\fsl}{{\mathfrak{sl}}}
\nc{\hsl}{{\widehat{\mathfrak{sl}}}}
\nc{\hgl}{{\widehat{\mathfrak{gl}}}}
\nc{\hg}{{\widehat{\mathfrak{g}}}}
\nc{\chg}{{\widehat{\mathfrak{g}}}{}^\vee}
\nc{\hn}{{\widehat{\mathfrak{n}}}}
\nc{\chn}{{\widehat{\mathfrak{n}}}{}^\vee}

\nc{\fA}{{\mathfrak{A}}}
\nc{\fB}{{\mathfrak{B}}}
\nc{\fD}{{\mathfrak{D}}}
\nc{\fE}{{\mathfrak{E}}}
\nc{\fF}{{\mathfrak{F}}}
\nc{\fG}{{\mathfrak{G}}}
\nc{\fK}{{\mathfrak{K}}}
\nc{\fL}{{\mathfrak{L}}}
\nc{\fM}{{\mathfrak{M}}}
\nc{\fN}{{\mathfrak{N}}}
\nc{\frP}{{\mathfrak{P}}}
\nc{\fS}{{\mathfrak S}}
\nc{\fU}{{\mathfrak{U}}}
\nc{\fZ}{{\mathfrak{Z}}}

\nc{\bb}{{\mathbf{b}}}
\nc{\bc}{{\mathbf{c}}}
\nc{\be}{{\mathbf{e}}}
\nc{\bj}{{\mathbf{j}}}
\nc{\bn}{{\mathbf{n}}}
\nc{\bp}{{\mathbf{p}}}
\nc{\bq}{{\mathbf{q}}}
\nc{\bfu}{{\mathbf{u}}}
\nc{\bv}{{\mathbf{v}}}
\nc{\bx}{{\mathbf{x}}}
\nc{\by}{{\mathbf{y}}}
\nc{\bw}{{\mathbf{w}}}
\nc{\bA}{{\mathbf{A}}}
\nc{\bB}{{\mathbf{B}}}
\nc{\bC}{{\mathbf{C}}}
\nc{\bK}{{\mathbf{K}}}
\nc{\bD}{{\mathbf{D}}}
\nc{\bH}{{\mathbf{H}}}
\nc{\bM}{{\mathbf{M}}}
\nc{\bN}{{\mathbf{N}}}
\nc{\bT}{{\mathbf{T}}}
\nc{\bV}{{\mathbf{V}}}
\nc{\bW}{{\mathbf{W}}}
\nc{\bX}{{\mathbf{X}}}
\nc{\bP}{{\mathbf{P}}}
\nc{\bZ}{{\mathbf{Z}}}

\nc{\sA}{{\mathsf{A}}}
\nc{\sB}{{\mathsf{B}}}
\nc{\sC}{{\mathsf{C}}}
\nc{\sD}{{\mathsf{D}}}
\nc{\sF}{{\mathsf{F}}}
\nc{\sK}{{\mathsf{K}}}
\nc{\sM}{{\mathsf{M}}}
\nc{\sO}{{\mathsf{O}}}
\nc{\sQ}{{\mathsf{Q}}}
\nc{\sP}{{\mathsf{P}}}
\nc{\sZ}{{\mathsf{Z}}}
\nc{\sfp}{{\mathsf{p}}}
\nc{\sr}{{\mathsf{r}}}
\nc{\sfb}{{\mathsf{b}}}
\nc{\sfc}{{\mathsf{c}}}
\nc{\sd}{{\mathsf{d}}}
\nc{\sfl}{{\mathsf{l}}}

\nc{\BK}{{\bar{K}}}

\nc{\tA}{{\widetilde{\mathbf{A}}}}
\nc{\tB}{{\widetilde{\mathcal{B}}}}
\nc{\tg}{{\widetilde{\mathfrak{g}}}}
\nc{\tG}{{\widetilde{G}}}
\nc{\TM}{{\widetilde{\mathbb{M}}}{}}
\nc{\tO}{{\widetilde{\mathsf{O}}}{}}
\nc{\tU}{{\widetilde{\mathfrak{U}}}{}}
\nc{\TZ}{{\tilde{Z}}}
\nc{\tx}{{\tilde{x}}}
\nc{\tbv}{{\tilde{\bv}}}
\nc{\tfP}{{\widetilde{\mathfrak{P}}}{}}
\nc{\tz}{{\tilde{\zeta}}}
\nc{\tmu}{{\tilde{\mu}}}

\nc{\urho}{\underline{\rho}}
\nc{\uB}{\underline{B}}
\nc{\uC}{{\underline{\mathbb{C}}}}
\nc{\ui}{\underline{i}}
\nc{\uj}{\underline{j}}
\nc{\ofP}{{\overline{\mathfrak{P}}}}
\nc{\oB}{{\overline{\mathcal{B}}}}
\nc{\og}{{\overline{\mathfrak{g}}}}
\nc{\oI}{{\overline{I}}}

\nc{\eps}{\varepsilon}
\nc{\hrho}{{\hat{\rho}}}

\nc{\one}{{\mathbf{1}}}
\nc{\two}{{\mathbf{t}}}

\nc{\Rep}{{\mathop{\operatorname{\rm Rep}}}}
\nc{\Tot}{{\mathop{\operatorname{\rm Tot}}}}
\nc{\Ker}{{\mathop{\operatorname{\rm Ker}}}}
\nc{\Hilb}{{\mathop{\operatorname{\rm Hilb}}}}
\nc{\End}{{\mathop{\operatorname{\rm End}}}}
\nc{\Ext}{{\mathop{\operatorname{\rm Ext}}}}
\nc{\RHom}{{\mathop{\operatorname{\rm RHom}}}}
\nc{\CHom}{{\mathop{\operatorname{{\mathcal{H}}\it om}}}}
\nc{\GL}{{\mathop{\operatorname{\rm GL}}}}
\nc{\gr}{{\mathop{\operatorname{\rm gr}}}}
\nc{\Id}{{\mathop{\operatorname{\rm Id}}}}
\nc{\defi}{{\mathop{\operatorname{\rm def}}}}
\nc{\length}{{\mathop{\operatorname{\rm length}}}}
\nc{\supp}{{\mathop{\operatorname{\rm supp}}}}

\nc{\Cliff}{{\mathsf{Cliff}}}
\nc{\Fl}{{\mathsf{Fl}}}
\nc{\Fib}{{\mathsf{Fib}}}
\nc{\Coh}{{\mathsf{Coh}}}
\nc{\FCoh}{{\mathsf{FCoh}}}

\nc{\reg}{{\text{\rm reg}}}

\nc{\cplus}{{\mathbf{C}_+}}
\nc{\cminus}{{\mathbf{C}_-}}
\nc{\cthree}{{\mathbf{C}_*}}
\nc{\Qbar}{{\bar{Q}}}

\nc{\bh}{{\bar{h}}}
\nc{\bOmega}{{\overline{\Omega}}}

\nc{\seq}[1]{\stackrel{#1}{\sim}}

\nc{\BBB}{{\mathbf B}}
\nc{\BBC}{{\mathbf C}}
\nc{\BBD}{{\mathbf D}}
\nc{\BBK}{{\mathbf K}}

\nc{\fV}{{\mathfrak{V}}}
\nc{\BG}{{\mathbb{G}}}
\nc{\oF}{{\overset{\circ}{\fF}}}

\theoremstyle{remark}
\newtheorem{Rem}{Remark}

\newcommand{\Qlbar}{{\overline {{\Bbb Q}_l} }}

\newcommand{\cal}{\mathcal}

\newcommand{\bu}{\bullet}

\newcommand{\To}{\longrightarrow}
\newcommand{\iso}{{\stackrel{\sim}{\longrightarrow}}}

\newcommand{\imbed}{\hookrightarrow}

\renewcommand{\P}{{\cal P}}

\renewcommand{\H}{{\cal H}}

\newcommand{\suml}{\sum\limits}
\newcommand{\oplusl}{\bigoplus\limits}

\renewcommand{\t}{{\frak t}}

\newcommand{\g}{{\frak g}}
\newcommand{\z}{{\frak z}}
\newcommand{\Lg}{\check{\mathfrak g}{}}
\newcommand{\Ll}{\check{\mathfrak l}{}}
\newcommand{\Ln}{\check{\mathfrak n}{}}
\newcommand{\Lp}{\check{\mathfrak p}{}}

\newcommand{\Lt}{\check{\mathfrak t}{}}
\newcommand{\Lb}{\check{\mathfrak b}{}}
\newcommand{\LN}{\check{N}{}}
\newcommand{\LG}{{\check{G}}}

\newcommand{\LL}{\check{L}{}}

\newcommand{\LT}{\check{T}{}}

\renewcommand{\O}{{\cal O}}

\newcommand{\cT}{{\cal T}}

\newcommand{\Zet}{{\Bbb Z}}
\newcommand{\Ce}{{\Bbb C}}

\newcommand{\Lotimes}{\overset{\rm L}{\otimes}}

\newcommand{\bI}{{\bf I}}
\newcommand{\GO}{{\bf G_O}}
\newcommand{\GK}{{\bf G_F}}
\newcommand{\GF}{{\bf G_F}}

\newcommand{\Aone}{{\mathbb A}^1}

\newcommand{\Gm}{{{\mathbb G}_m}}

\newcommand{\Fqbar}{{\bar{\mathbb F}_q}}

\newcommand{\Sa}{S}
\newcommand{\Satil}{{\tilde S}}

\newcommand{\HCh}{\HC_\hbar}
\newcommand{\HChfr}{\HC_\hbar^{fr}}

\newcommand{\kappah}{{\kappa_\hbar}}

\newcommand{\kk}{\mathsf k}

\newcommand{\beq}{\begin{equation}}
\newcommand{\eeq}{\end{equation}}

\renewcommand{\proof}{{\it Proof }}

\newcommand{\eqss}{{\mathcal{S}}}
\newcommand{\fufu}{{\mathfrak{F}}}
 \newcommand{\CatTopss}{{\mathcal{IC}}}
  \newcommand{\Bimodfr}{{\HC_\hbar^{fr}}}
  \newcommand{\qkap}{{\kappa_\hbar}}
  \newcommand{\ds}{{\mathfrak D}}
 \newcommand{\Ver}{{\mathbb D}}

 \newcommand{\prd}{{\mathbf D}}
 \newcommand{\inv}{{\mathfrak C}}

\newcommand{\psf}{{\kappa_\hbar}}
\newcommand{\cK}{{\mathcal K }}
\newcommand{\bimg}{{\ Bimodgr\  }}
\newcommand{\DMS}{{\ Satak-model\  }}

\newcommand{\PerSatak}{{Perv_{\GO}(\Gr)}}
\newcommand{\DerSatak}{{D_{\GO\rtimes\Gm}(\Gr)}}
\newcommand{\DerSatakneloop}{{D_{\GO}(\Gr)}}

\newcommand{\PerSatake}{{Perv_{\GO\rtimes\Gm}(\Gr)}}

\newcommand{\Dsym}{{\Sym^{[]}(\Lg)}}
\newcommand{\Duh}{{U_\hbar^{[]}}}
\newcommand{\la}{\lambda}

\newcommand{\Epur}{{\mathcal E}}

\author{Roman Bezrukavnikov and Michael Finkelberg}
\title{Equivariant Satake category and Kostant-Whittaker reduction}

\thanks{{\em 2000 Mathematics Subject Classification} 19E08 (primary),
  22E65, 37K10 (secondary).\\
{\em Keywords:} affine Grassmannian, Langlands dual group, Toda lattice.}
\dedicatory{To Victor Ginzburg on his birthday}
\begin{document}
\maketitle
\tableofcontents

\begin{abstract}
We explain (following V.~Drinfeld)
how the $G(\Ce[[t]])$ equivariant derived category
of the affine Grassmannian can be described in terms
of coherent sheaves on the Langlands dual Lie algebra equivariant
with respect to the adjoint action, due to some old results of V.~Ginzburg.
The global cohomology functor corresponds under this identification to restriction
to the Kostant slice.
We extend this description to loop rotation equivariant
derived category, linking it to Harish-Chandra bimodules
for the Langlands dual Lie algebra, so that the global cohomology functor corresponds to
the quantum Kostant-Whittaker reduction of a Harish-Chandra bimodule.
We derive a conjecture of \cite{BFM}, which identifies
the loop-rotation equivariant homology of the affine Grassmannian
with quantized Toda lattice.
\end{abstract}

\section{Introduction}
Let $G$ be a semi-simple algebraic group over an algebraically closed 
characteristic zero field $\kk$. The fundamental object of the 
geometric Langlands duality theory is the so-called {\em Satake 
category $\PerSatak$}. The latter is defined as the category of perverse
sheaves on the affine (loop) Grassmannian $\Gr$
equivariant with respect to the group of regular loops $\GO$.

It turns out that convolution provides $\PerSatak$ with a tensor structure,
and the celebrated geometric Satake isomorphism theorem establishes an
equivalence between 
$\PerSatak$ and the category of representations
of the Langlands dual group $\LG$.

By its very definition $\PerSatak$ arises as the heart of the t-structure
on a monoidal triangulated category -- the equivariant derived category
$\DerSatakneloop$. It is a natural question (raised, in particular,
 by V.~Drinfeld) 
to describe $\DerSatakneloop$
 in terms of the dual group. Drinfeld has also noticed
that at least some form of the answer\footnote{We do not address 
Drinfeld's problem to find a more natural derivation of the description,
making compatibility with finer structures transparent. We understand
that D.~Gaitsgory and J.~Lurie have made a significant progress in this
 direction.}
follows from the results of V.~Ginzburg's preprint~\cite{G}. 
In the present paper we reproduce this description and extend
it to a description of the {\em loop rotation equivariant} derived Satake
category $\DerSatak$.

The description of $\DerSatakneloop$ links it to 
{\em conjugation equivariant coherent sheaves} on the Langlands dual Lie
algebra. The additional $S^1$
(or $\Gm$) equivariance is connected to quantization
of these to Harish-Chandra bimodules (see Theorem~\ref{bone} for
a precise formulation).

The argument follows the strategy of~\cite{G}; 
it is based on another result of Ginzburg~\cite{G0},
which reduces the question to computation of the global equivariant
cohomology of IC sheaves as modules over the global equivariant
cohomology algebra $H^\bullet_{\GO\rtimes\Gm}(\Gr)$.
By an explicit calculation we show that $H^\bullet_{\GO\rtimes\Gm}(\Gr)$ 
is related to the 
tensor square of the center of the enveloping of the dual Lie algebra
$\Lg$, while the global cohomology modules correspond to the bimodules,
which describe twisting with a finite dimensional $\LG$-modules
on the category of {\em Whittaker modules}.
%are identified with
%{\em Kostant-Whittaker reduction}.
 This allows us to  relate $\DerSatak$
to Harish-Chandra bimodules, so that the global cohomology is identified
with the Kostant-Whittaker reduction.

As an application we prove a conjecture of~\cite{BFM} which
identifies the algebra of global equivariant homology
of $\Gr$ equipped with the convolution algebra structure with
the quantized Toda lattice (Theorem~\ref{TodaTh}).
Note that the quantized Toda lattice also appears in the apparently related
computations by Givental, Kim and others of quantum $D$-module 
(quantum cohomology) of the flag variety of $G$, see e.g.~\cite{Kim}.

D.~Ben-Zvi and D.~Nadler have informed us that they have a more conceptual
proof of some of our results, see~\cite{BN}.

{\bf Acknowledgements.} As should be clear from the above, 
this modest token of gratitude to Victor
Ginzburg is largely an outgrowth of his own works.
Throughout the mathematical biography of the authors
his ingenuity, generosity and enthusiasm have been an important source of
 support and inspiration.
We wish Vitya many happy returns of the day and a lot of exciting
mathematics ahead!

This paper was partly inspired by discussions 
with V.~Drinfeld around 2001, and some of the techniques were 
explained to us by L.~Positselski around 1998.
We also benefited a lot from explanations by
A.~Beilinson, D.~Gaitsgory, J.~Lurie and V.~Lafforgue. We gratefully
acknowledge all these contributions.
 We also thank D.~Ben-Zvi and D.~Nadler for their interest
in this work and for sharing their related preprint~\cite{BN}.

Work on this project has been partially supported by a
DARPA grant HR0011-04-1-0031.
 R.B. was partially supported by Sloan Foundation and NSF grant 
DMS-0625234,
and visited Independent University of Moscow during the work;
M.F. was partially supported by the CRDF award RUM1-2694, and
the ANR program ``GIMP'', contract number ANR-05-BLAN-0029-01,
and visited MSRI during the work.

%American Embassy. 
%nevinno ubiennye irakcy i bezvremenno ohrenevshie ostal'nye.

\section{Notation and statements of the results}

\subsection{Notation}

$\kk$ is an algebraically closed characteristic zero coefficient
field. Let $G$ be a semisimple complex algebraic group, 
$\GO=G(\BC[[t]]),\ \GK=G(\BC((t)))$. 
The affine Grassmannian $\Gr=\Gr_G=\GK/\GO$
carries the category $Perv_\GO(\Gr)$ of $\GO$-equivariant perverse 
constructible sheaves. It is equipped with the convolution monoidal 
structure, and is tensor equivalent to the tensor category $Rep(\LG)$ of 
representations of the Langlands dual group $\LG$ over the field $\kk$
(see~\cite{L},~\cite{BD},~\cite{G}). We denote by $\tilde{S}:\ 
Rep(\LG)\to Perv_\GO(\Gr)$ the geometric Satake isomorphism functor,
and by $S:\ Rep(\LG)\to Perv_{\GO\rtimes\Gm}(\Gr)$ its extension 
to the monoidal category of $\GO\rtimes\Gm$-equivariant perverse
constructible sheaves.
The Lie algebra of $\LG$ is denoted by $\Lg$. We choose a Cartan subalgebra
$\Lt\subset\Lg$; the corresponding Cartan torus in $\LG$ is denoted $\LT$.
We choose the opposite Borel subalgebras $\Lb_\pm\supset\Lt$ with 
nilpotent radicals $\Ln_\pm$ and corresponding unipotent subgroups 
$\LN_\pm\subset\LG$. The Weyl group of $\LG$ is denoted by $W$.

Let $e,h,f\in \Lg$ be a principal ${\mathfrak{sl}}_2$ triple
such that $f\in\Ln_-,\ e\in\Ln_+$. We have the
Kostant slice $e+\z(f)$ to the
principal nilpotent orbit. It is known that $e+\z(f)\iso \Lg/Ad (\LG)$,
and also $e+\z(f)\iso (e+\Lb_-)/\LN_-$; moreover the $\LN_-$ action on
$e+\Lb_-$ is free. 
Let $\Sigma$, $\Upsilon$ be the images of $e+\z(f)$, $e+\Lb_-$
under a $\LG$-invariant isomorphism $\Lg\cong \Lg^*$. Thus we have
$\Sigma\iso \Upsilon/\LN_-=\Lt^*/W$ canonically.

The total space of the tangent bundle of $\Lt^*/W$ is denoted $\bT(\Lt^*/W)$.

\subsection{Asymptotic $\Lg$-modules}
\label{agm}
Let $U=U(\Lg)$ be the enveloping algebra, and let  $U_\hbar$ be the
``graded enveloping" algebra, i.e. the graded $\kk[\hbar]$-algebra
generated by $\g$ with relations $xy-yx = \hbar [x,y]$ for $x,y
\in \Lg$ (thus $U_\hbar$ is obtained from $U$ by the standard Rees
construction which produces a graded algebra from a filtered one);
the adjoint action extends to the action of $\LG$ on $U_\hbar$. We
define the category $\widetilde \HCh$ of ``$\hbar$-Harish-Chandra
 bimodules" as follows: an object $M$ of this category is a 
graded $U_\hbar^2:=U_\hbar
 \otimes_{\kk[\hbar]} U_\hbar\simeq U_\hbar\otimes_\kk U$-module 
equipped with an algebraic
 action $\rho$ of $\LG$ such that: (1) the action map
 $U_\hbar\otimes_{\kk[\hbar]}U_\hbar \otimes M \to M$
 is $\LG$ equivariant, and (2) for $x\in \Lg$ the action of
 $(x\otimes 1 +1\otimes x)\in U_\hbar \otimes _{\kk[\hbar]} U_\hbar$
coincides with $\hbar \cdot d\rho(x)$. 
The functor of restriction
from $U_\hbar \otimes U_\hbar$ to $U_\hbar\otimes 1$ is an
equivalence between $\widetilde \HCh$ and the category of $\LG$-modules
equipped with an equivariant $U_\hbar$-action; the same is true
for the restriction to $1\otimes U_\hbar$. We let $\HCh \subset
\widetilde \HCh$ denote the full subcategory of objects which are
finitely generated as $U_\hbar\otimes 1$ modules (equivalently, as
$1\otimes U_\hbar$ modules).

Notice that the full
subcategory of $\HCh$ consisting of objects where $\hbar$ acts by
zero is identified with the category $Coh^\LG(\Lg^*)$ of coherent
sheaves on $\Lg^*$ equivariant under the coadjoint action; while
for $s\in \kk$, $s\ne 0$ the subcategory where $\hbar$ acts by $s$
is identified with the category of Harish-Chandra bimodules.

%\begin{Rem}
%\label{mod-bimod}
%For $M\in \HCh$ the action of $U_\hbar \otimes U_\hbar$ and the 
%differential of the $\LG$ action together generate an action of
%the algebra $U_\hbar^2$ isomorphic to $U_\hbar \otimes _\kk U$
%(the smash product with respect to the adjoint action of $U$ on 
%$U_\hbar$). Thus $M$ can be viewed as a
%$\LG$-equivariant $U_\hbar$-module.
%\end{Rem}

\subsection{Kostant functor $\kappah$} 
\label{Kostant}
We now proceed to
define a functor $\kappah:\HCh\to QCoh^\Gm((\Lt^*/W)^2\times \Aone)$.

Let $\psi:U_\hbar (\Ln_-)\to\kk[\hbar]$ be a homomorphism
such that $\psi(f_\alpha)=1$ for any simple root 
$\alpha$, and a root generator $f_\alpha\in\Ln_-\subset U_\hbar(\Ln_-)$.

Define  $U_\hbar^2(\Ln_-)\subset U_\hbar ^2$ by 
$U_\hbar^2(\Ln_-)=U_\hbar(\Ln_-)
\otimes U(\Ln_-)$. We extend $\psi$ to a  
character $\psi_{(2)}: U_\hbar^2(\Ln_-)
=U_\hbar(\Ln_-) \otimes U(\Ln_-)
\to\kk[\hbar]$ trivial on the second multiple. Note that its restriction
to the first copy of $U_\hbar$ is $\psi$, and its restriction
to the second copy is $(-\psi)$. 

%We set $\kappahbar(M)=M\Lotimes_{U_\hbar^2(\Ln)} \psi_{(2)}[\dim \Ln]$.

We set $\kappah(M)=(M\Lotimes_{U_\hbar(\Ln_-)_2}(-\psi))^{\LN_-}$
where the action of the second copy of $U_\hbar$ is used (though using 
the first one we get a canonically isomorphic functor). Clearly,
$\kappah(M)$ is equipped with the action of the Harish-Chandra center
$Z(U_\hbar)\otimes_{\kk[\hbar]}Z(U_\hbar)=
\CO((\Lt^*/W)\times(\Lt^*/W)\times\Aone)$, 
and with the grading (coming from the action of the Cartan element $h$
of the principal $\mathfrak{sl}_2$), so we may
view $\kappah(M)$ as a $\Gm$-equivariant quasicoherent sheaf on 
$(\Lt^*/W)^2\times\Aone$.

If $X$ is a scheme, and $Z\subset X$ is a closed subscheme let
$N_XZ$ be the {\it deformation to the normal cone} to $Z$, see
\cite{Fu}. It is equipped with a morphism $N_XZ\to X\times\Aone$
(with coordinate $\hbar$ on $\Aone$),
and is defined as the relative spectrum of the sheaf of subalgebras in
$\CO_X[\hbar^{\pm1}]$ generated by the elements $f\hbar^{-1},\ 
f\in\CO_X:\ f|_Z=0$. In other words, $N_XZ$ is the affine blowup
$\operatorname{Bl}_{Z\times\{0\}}^{\operatorname{aff}}(X\times\BA^1)$ of $X\times\BA^1$ at $Z\times\{0\}$:
the complement of the strict transform of $X\times\{0\}$ in the blowup $\operatorname{Bl}_{Z\times\{0\}}(X\times\BA^1)$.

If $M$ is a free $\kk[\hbar]$-module, 
the action of $\CO((\Lt^*/W)^2\times\Aone)$
on $\kappah(M)$ extends uniquely to the action of
$\CO(N_{(\Lt^*/W)^2}\Delta)$ where $\Delta\subset(\Lt^*/W)^2$ is the 
diagonal. So we can and will view $\kappah(M)$ as a $\Gm$-equivariant
coherent sheaf on $N_{(\Lt^*/W)^2}\Delta$.

For $V\in Rep(\LG)$ we define the $\hbar$-Harish-Chandra bimodule
$Fr(V)$ by $Fr(V)=U_\hbar\otimes V$ with its natural $\LG$-module
structure $g(y\otimes v)=Ad(g)(y)\otimes g(v)$, and with the
$U_\hbar \otimes U_\hbar$ action specified by $x\otimes u(y\otimes
v)= xyu\otimes v +\hbar\cdot xy \otimes u(v)$ for $x,u\in \Lg
\subset U_\hbar$. In other words, $Fr(V)$ is obtained by applying the induction
(left adjoint to the restriction functor) $Rep(\LG)\to \HCh$ to $V$.
We set $\phi(V):=\kappa_\hbar(Fr(V))$.

 Clearly  $Fr(V)$ is a projective object of
$\HCh$ for any $V\in Rep(\LG)$; we call an object of the form
$Fr(V)$ a free $\hbar$-Harish-Chandra bimodule. We define the full
subcategory  $\HChfr\subset \HCh$ to consist of all free objects.

\subsection{Equivariant cohomology of $\Gr_G$}

Note that $H^\bullet_{\GO\rtimes\Gm}(\Gr_G)=
H^\bullet_{\Gm}(\GO\backslash\GF/\GO)$,
whence two morphisms $pr_1^*, pr_2^*:\ \O(\t/W)=H^\bullet_\GO(pt)\to
H^\bullet_{\GO \rtimes \Gm}(\Gr_G)$.
We also have a morphism $pr^*:\ \kk[\hbar]=H^\bullet_\Gm(pt)\to
H^\bullet_{\GO \rtimes \Gm}(\Gr_G)$.

\begin{thm}
\label{normocon}
a) Assume $G$ is simply connected. 
We have a canonical isomorphism $H^\bullet_{\GO \rtimes \Gm}(\Gr_G)\cong
\O(N_{\Lt^*/W\times \Lt^*/W}\Delta)$ where $\Delta \subset (\Lt^*/W)^2$ is
the diagonal. Here the projection $N_{ (\Lt^*/W)^2}\Delta\to \Aone$
corresponds to the homomorphism $H_{\Gm}^\bu(pt)\to H^\bu_{\GO
\rtimes \Gm} (\Gr_G)$; and the two projections $N_{ (\Lt^*/W)^2}\Delta
\to \Lt^*/W=\t/W$ correspond to the two homomorphisms
$H_{\GO}^\bu(pt)\to H^\bu_{\GO
\rtimes \Gm} (\Gr_G)$. The isomorphism is specified uniquely
 by these requirements.

b) For arbitrary $G$ we have a canonical isomorphism
$H^*_{\GO \rtimes \Gm}(\Gr_G)\cong\bigoplus_{\pi_1(G)}\O(N_{\Lt^*/W\times \Lt^*/W}\Delta)$.
\end{thm}

\begin{Rem}
To simplify the exposition we assume from now on that $G$ is simply
connected. 
\end{Rem}

Cohomology of any complex of sheaves on a topological space
carries an action of the cohomology algebra of the space; thus we have
the functor of equivariant cohomology
$$H_{\GO\rtimes\Gm}^\bullet: D_{\GO\rtimes\Gm}(\Gr)\to 
H^\bullet_{\GO\rtimes\Gm}(\Gr)
-mod^{gr} = Coh^{\Gm}(N_{(\Lt^*/W)^2}\Delta)$$ 
where $D_{\GO\rtimes\Gm}(\Gr)$ denotes the bounded constructible
equivariant derived category, and the grading on 
$H^*_{\GO\rtimes\Gm}(\Gr,?)$ is the one by the cohomology degree. 

\begin{thm}
\label{free_cohh}
The functor $\Sa: Rep(\LG)\to Perv_{\GO\rtimes \Gm}(\Gr)$ 
extends to a full imbedding $\Sa_\hbar:\HChfr \to
D_{\GO\rtimes \Gm}(\Gr)$, such that 
\begin{equation}\label{cogokaph}
\kappah \cong H_{\GO\rtimes\Gm}^\bullet\circ \Sa_\hbar.
\end{equation}
Such an extension $\Sa_\hbar$ (for a  fixed isomorphism \eqref{cogokaph})
is unique.
\end{thm}

\subsection{Equivariant homology and quantum Toda lattice}
Let $\CalD_\hbar(\LG)$ stand for the sheaf of $\hbar$-differential
operators on $\LG$: its global sections is the smash product of 
$U_\hbar$ and $\CO(\LG)$. The action of $\Ln_-$ by the left-invariant
(resp. right-invariant) vector fields on $\LG$ gives rise to the
homomorphism $\sfl$ (resp. $\sr$): $U_\hbar(\Ln_-)\to\CalD_\hbar(\LG)$.
Let $I_\psi\subset\CalD_\hbar(\LG)$ be the left ideal generated by the
$\hbar$-differential operators of the sort
$\sfl(u_1)-\psi(u_1)+\sr(u_2)+\psi(u_2);\ u_1,u_2\in U_\hbar(\Ln_-)$.
We consider the quantum hamiltonian reduction 
$$(\CalD_\hbar(\LG)/I_\psi)^{\LN_-\times\LN_-}$$
where the first (resp. second) copy of $\LN_-$ acts on $\LG$ (and hence
on $\CalD_\hbar(\LG)$) by the left
(resp. right) translations: $(n_1,n_2)\circ g:=n_1gn_2^{-1}$. 
It is an algebra containing a commutative subalgebra $Z(U_\hbar)$
(via the embedding $Z(U_\hbar)\hookrightarrow\CalD_\hbar(\LG)$ as both
left- and right-invariant $\hbar$-differential operators).
Note that the action of $\LN_-\times\LN_-$ on the ``big Bruhat cell'' 
$C_{w_0}:=\LN_-\cdot\LT\cdot w_0\cdot\LN_-\subset\LG$ is free, and
hence the quantum hamiltonian reduction of $\CalD_\hbar(C_{w_0})$ is
isomorphic to $\CalD_\hbar(\LT)$. This is the classical
Kazhdan-Kostant construction of the quantum Toda lattice,
see~\cite{K2}. Thus, the quantum Toda lattice is a certain
localization of $(\CalD_\hbar(\LG)/I_\psi)^{\LN_-\times\LN_-}$.
In what follows we will call
$(\CalD_\hbar(\LG)/I_\psi)^{\LN_-\times\LN_-}$ ``the quantized Toda
lattice'', somewhat abusing the language. 
The following result was conjectured in~\cite{BFM}.

\begin{thm}
\label{TodaTh}
The convolution algebra of equivariant homology
$H^{\GO\rtimes \Gm}_\bullet(\Gr)$ is naturally isomorphic to the
quantized Toda lattice $(\CalD_\hbar(\LG)/I_\psi)^{\LN_-\times\LN_-}$.
The embedding $Z(U_\hbar)\simeq H_{\GO\rtimes\Gm}(pt)\hookrightarrow 
H^{\GO\rtimes\Gm}_\bullet(\Gr)$ corresponds to the embedding
$Z(U_\hbar)\hookrightarrow(\CalD_\hbar(\LG)/I_\psi)^{\LN_-\times\LN_-}$.
\end{thm}

\subsection{Quasiclassical limit}
\label{quasic}
Recall that the fiber of $N_XZ$ over $0\in\Aone$ is the normal cone to $Z$
in $X$. In particular, the fiber of $N_{(\Lt^*/W)^2}\Delta$ over
$0\in\Aone$ is the total space of the tangent bundle $\bT(\Lt^*/W)$.
Thus, Theorem~\ref{normocon} implies the canonical isomorphism
$H^\bullet_\GO(\Gr)\simeq\CO(\bT(\Lt^*/W))$. On the other hand, 
$H^\bullet_\GO(\Gr)$ was computed by V.~Ginzburg in~\cite{G} in terms of
the universal centralizer bundle of $\Lg$. The two computations are related
as follows. 

The variety $(\Lg^*)^{reg} $ of regular elements in $\Lg^*$
carries a sheaf of commutative Lie algebras $\z \subset \Lg\otimes
\O$ whose fiber at a point $\xi\in (\Lg^*)^{reg}$ is the
stabilizer of $\xi$. We claim a canonical isomorphism $\z\cong
pr^*(\cT^*)$ where $pr:(\Lg^*)^{reg} \to \Lt^*/W$ is the
projection to the spectrum of invariant polynomials, and $\cT^*$
stands for the cotangent sheaf. Indeed, the fiber of $pr^*(\cT^*)$
at a point $\xi\in \Lg^*$ is dual to  the cokernel of the map $\Lg \to
\Lg^*$, $x\mapsto coad(x)(\xi)$; thus it is canonically isomorphic
to the kernel of the dual map (which happens to coincide with the
original map), which is exactly the fiber of $\z$ at $\xi$.

In view of this identification, one should compare Lemma~\ref{MIt}
in subsection~\ref{secr} below with Ginzburg's description in~\cite{G}
of $\GO$-equivariant Intersection Cohomology
of a $\GO$-orbit in $\Gr$ as a $\z$-module.

We now proceed to
define the {\em Kostant functor} $\kappa: Coh^{\LG\times\Gm}(\Lg^*)\to Coh^\Gm
(\bT(\Lt^*/W))$, $\Gm$-equivariant coherent sheaves on the tangent bundle to
$\Lt^*/W$.

 If $\F \in Coh^\LG(\Lg^*)$ is equipped with
an equivariant structure, then $\F|_{(\Lg^*)^{reg} }$ carries an
action of $\z$; thus by the previous paragraph it defines a
coherent sheaf on the total space of the pull-back of the tangent
bundle under $pr$. Restricting this sheaf to the preimage of
$\Sigma$ we get a coherent sheaf on the tangent bundle to
$\Sigma=\Lt^*/W$ which we denote by $\overline{\kappa}(\F)$. Notice that
$\overline{\kappa}(\F)= (\F|_{\Upsilon})^{\LN_-}$ (where we do not distinguish
between a coherent sheaf on an affine variety and the module of
its global sections). An obvious modification of this definition
yields a functor $\kappa: Coh^{\LG\times \Gm}(\Lg^*)\to Coh^{\Gm}
(\bT(\Lt^*/W))$ (where the action of $\Gm$ on $\Lt^*/W$ is the natural one).

Define the full subcategory  
%$Coh_{fr}^{\LG}(\Lg^*)\subset Coh^{\LG}(\Lg^*)$,
 $Coh_{fr}^{\LG\times \Gm}(\Lg^*)\subset Coh^{\LG\times \Gm}(\Lg^*)$ 
to consist of all objects of the form 
$V\otimes \O_{\Lg^*}$, for $V\in Rep(\LG\times \Gm)$.

Recall that $\Satil:\ Rep(\LG)\to Perv_{\GO}(\Gr)$ 
is the composition of $\Sa$
with the forgetful functor $Perv_{\GO\rtimes \Gm}(\Gr)\to Perv_{\GO}(\Gr)$.
We have the functor of $\GO$-equivariant cohomology
$$H_\GO^\bullet: D_{\GO}(\Gr)\to H^\bullet_{\GO}(\Gr)
-mod^{gr} = Coh^{\Gm}(\bT(\Lt^*/W))$$

\begin{thm}\label{free_coh} The functor $\Satil$
extends to a full imbedding $\Satil_{qc}:\ Coh_{fr}^{\LG\times \Gm}(\Lg^*)\to
D_{\GO}(\Gr)$, such that there exists an isomorphism 
\begin{equation}\label{cogokap} 
\kappa \cong H_\GO^\bullet\circ \Satil_{qc}.
\end{equation}
Such an extension $\Satil_{qc}$ (for a  fixed isomorphism \eqref{cogokap})
is unique.
\end{thm}

\subsection{Equivalences}
To  a differential graded algebra $A$ one can associate
the triangulated  category $D(A)$  of differential graded modules localized
by quasi-isomorphism; and a full triangulated subcategory 
$D_{perf}(A)\subset D(A)$ 
of perfect complexes. Thus $D_{perf}(A)$ is
the full subcategory in the latter category consisting of perfect complexes
(i.e. generated by the free module under cones and direct summands).
Given an algebraic group $H$ acting
on a dg-algebra $A$, we can consider equivariant dg-modules and localize
them by quasi-isomorphisms, arriving at the equivariant version
$D^H_{perf}(A)$.

 We now consider the ``differential-graded versions'' $\Dsym$, $\Duh$
of the graded algebras 
$\Sym (\Lg)$, $U_\hbar(\Lg)$. By definition $\Dsym$, $\Duh$
are differential graded algebras with zero differential, which as algebras are
isomorphic to $\Sym(\Lg)$, $U_\hbar(\Lg)$ respectively. The cohomological
grading is defined so that elements of $\Lg$ and $\hbar$ have degree two.
Recall from section \ref{agm} that an
asymptotic Harish-Chandra bimodule $M\in\HC_\hbar$ is nothing but a
$\LG$-equivariant $\Duh$-module. Using this identification we can
transfer tensor product of asymptotic Harish-Chandra bimodules to
a monoidal structure on the category of $\LG$-equivariant
$\Duh$-modules. It gives rise to a monoidal structure on
$D^{\LG}_{perf}(\Duh)$. Similarly, we define a monoidal structure on 
$D^{\LG}_{perf}(\Dsym)$.

\begin{thm}\label{polkovnik}
There exist canonical equivalences of monoidal triangulated categories
$D^{\LG}_{perf}(\Duh) \cong D_{\GO\ltimes \Gm}(\Gr)$,
$D^{\LG}_{perf}(\Dsym) \cong D_{\GO}(\Gr)$.

\end{thm}

The following statement is an immediate consequence of the Theorem,
which has the (psychological) advantage  of bypassing 
the notion of a dg-algebra.

\begin{cor}
\label{bone}
a) The derived
graded Harish-Chandra bimodule category is
Koszul dual to a graded version of the 
loop-rotation equivariant Satake category,
i.e. (cf.~\cite{BGS},~\cite{ABG},~\cite{Hu}) there exists a functor 
$\Phi:D^b(\HC_\hbar) \to D_{\GO\rtimes \Gm}(\Gr)$, such that 

i) $\Phi(M(1)) \cong \Phi(M)[1]$ in a natural way.

ii) For $M_1,\, M_2 \in D^b(\HC_\hbar)$, $\Phi$ induces an isomorphism
$$\suml_{n,m}Hom(M_1, M_2(n)[m])\cong \suml_k Hom(\Phi(M_1),\Phi(M_2)[k])$$

iii)
The image of $\Phi$ generates the target category as a triangulated 
category.

Moreover, $\Phi$ carries a natural monoidal structure, and  
$\Phi|_{\HC_\hbar^{fr}}$ coincides with $S_\hbar$.

b) The derived category of graded $\LG$-equivariant $\Sym(\Lg)$
modules is 
Koszul dual to a graded version of the  Satake category,
i.e. there exists a functor 
$\Phi':D^b(Coh^{\LG\times \Gm}(\Lg) ) \to D_{\GO}(\Gr)$, 
satisfying the properties similar to those listed in part (a).
This functor has a natural monoidal structure.
\end{cor}

\section{Topology}

\subsection{Proof of Theorem~\ref{normocon}}
a) First we construct the morphism $\alpha:\
\O(N_{\Lt^*/W\times \Lt^*/W}\Delta)\to H^\bullet_{\GO \rtimes \Gm}(\Gr)$.
Recall that $H^\bullet_{\GO\rtimes\Gm}(\Gr)=
H^\bullet_{\Gm}(\GO\backslash\GF/\GO)$,
whence two morphisms $pr_1^*, pr_2^*:\ \O(\t/W)=H^\bullet_\GO(pt)\to
H^\bullet_{\GO \rtimes \Gm}(\Gr)$.
We also have a morphism $pr^*:\ \kk[\hbar]=H^\bullet_\Gm(pt)\to
H^\bullet_{\GO \rtimes \Gm}(\Gr)$. Since 
$H^\bullet_{\GO \rtimes \Gm}(\Gr)|_{\hbar=0}=
H^\bullet_\GO(\Gr)$, it follows that $pr_1^*|_{\hbar=0}=pr_2^*|_{\hbar=0}$.
Hence the morphism $(pr_1^*,pr_2^*,pr^*):\ \O(\t/W\times\t/W\times\Aone)\to
H^\bullet_{\GO \rtimes \Gm}(\Gr)$ factors through the desired morphism
$\O(\t/W\times\t/W\times\Aone)\to
\O(N_{\t/W\times \t/W}\Delta)\overset{\alpha}{\to} 
H^\bullet_{\GO \rtimes \Gm}(\Gr)$.

Next we prove that $\alpha$ is an embedding. It suffices to prove that
the localized morphism $$\alpha_{loc}:\ 
\O(N_{\t/W\times \t/W}\Delta)\otimes_{\O(\t/W\times\Aone)}
\operatorname{Frac}(\O(\t\times\Aone))
\to H^\bullet_{\GO \rtimes \Gm}(\Gr)\otimes_{\O(\t/W\times\Aone)}
\operatorname{Frac}(\O(\t\times\Aone))$$ is an embedding. 
Note that the RHS is the 
localized equivariant cohomology $H^\bullet_{T\times\Gm}(\Gr)_{loc}$,
which embeds into the inverse limit of the localized equivariant cohomology
$H^\bullet_{T\times\Gm}(\Gr_\lambda)_{loc}$ of the $\GO$-orbit
closures $\Gr_\lambda\subset\Gr$.
By the Localization Theorem for torus-equivariant cohomology,
the latter is $\prod_\lambda\operatorname{Frac}(\O(\t\times\Aone))$
(cf.~\ref{can}), and $\alpha_{loc}$ is an embedding.

Finally, it remains to check that the graded dimensions of
$H^\bullet_{\GO \rtimes \Gm}(\Gr)$, and of
$\O(N_{\t/W\times \t/W}\Delta)$ coincide. Here the grading of
$\O(N_{\t/W\times \t/W}\Delta)$ comes from the natural $\Gm$-actions on
$\t$ and $\Aone$. To this end note that the graded dimension of 
$H^\bullet_{\GO \rtimes \Gm}(\Gr)$ coincides with that of
$H^\bullet_\GO(pt)\otimes H^\bullet_\Gm(pt)\otimes H^\bullet(\Gr)$,
that is $\kk[x_1,\ldots,x_r,y_1,\ldots,y_r,\hbar]$. Here $r$ is the rank of
$G$; the degree of $\hbar$ is 2; the degrees of $x$'s are twice the 
exponents of $\Lg$; the degrees of $y$'s are twice the exponents of
$\Lg$ minus 2.

Now to compute the graded dimension of $\O(N_{\t/W\times \t/W}\Delta)$ we
use that $\t/W$ is isomorphic to the vector space $\Sigma$. More generally,
for vector spaces $V,V'$ we have an isomorphism 
$\beta: V\times V'\times\Aone\iso N_{V\times V'}V$. In effect, the map
$\gamma: V\times V'\times\Aone\to V\times V'\times\Aone,\
(v,v',a)\mapsto(v,av',a)$, factors through the desired isomorphism:
$V\times V'\times\Aone\overset{\beta}{\to}N_{V\times V'}V\to 
V\times V'\times\Aone$. We derive an isomorphism 
$N_{\t/W\times\t/W}\Delta\cong\t/W\times\t/W\times\Aone$
which lowers the $\Gm$-weights in the second copy of $\t/W$ by 1, 
whence the desired formula for the graded
dimension of $\O(N_{\t/W\times\t/W}\Delta)$.

This completes the proof of the part a) of the theorem.

b) For general $G$ (not necessarily simply connected), we will describe $H^\bullet_{\GO\rtimes\Gm}(\Gr)$ as $W$-invariants in $H^\bullet_{T\times\Gm}(\Gr_G)$.
First we describe the equivariant cohomology ring $H^\bullet_{T\times\Gm}(\Gr_{G^{\on{sc}}})$ of the neutral
connected component of $\Gr_G$. To this end we consider the affinized Cartan subalgebra $\t_{\on{aff}}=\t\times\BA^1$
(with extra coordinate $\hbar$) with its natural projection to $\BA^1$. Let
$\Gamma\subset\t_{\on{aff}}\times_{\BA^1}(\t_{\on{aff}}/W)$ be the graph of the projection
$\t_{\on{aff}}\to\t_{\on{aff}}/W$ (note that the action of $W$ is trivial on $\BA^1$). The argument in the proof
of~a) shows that
$H^\bullet_{T\times\Gm}(\Gr_{G^{\on{sc}}})=\O\left(\on{Bl}^{\on{aff}}_\Gamma(\t_{\on{aff}}\times_{\BA^1}(\t_{\on{aff}}/W))\right)$.
Note that the RHS has a natural $W$-action. Taking $W$-invariants, we reproduce the answer in~a):
$H^\bullet_{{\mathbf G}^{\on{sc}}_{\mathbf O} \rtimes \Gm}(\Gr_{G^{\on{sc}}})\cong
\Big(\O\left(\on{Bl}^{\on{aff}}_\Gamma(\t_{\on{aff}}\times_{\BA^1}(\t_{\on{aff}}/W))\right)\Big)^W\cong\O(N_{\Lt^*/W\times \Lt^*/W}\Delta)$.

For general $G$ (not necessarily simply connected), the coweight lattice $\Lambda=X_*(T)$ contains the coroot
lattice $Q$. We have a natural action of $\Lambda$ on $\t_{\on{aff}}$ : $\lambda(x,c\hbar)=(x+c\lambda,c\hbar)$,
compatible with the action of $W$ on $\Lambda$ and on $\t_{\on{aff}}$. Thus we obtain an action of $\Lambda$
on $\O\left(\on{Bl}^{\on{aff}}_\Gamma(\t_{\on{aff}}\times_{\BA^1}(\t_{\on{aff}}/W))\right)$. In particular, the coroot lattice
$Q\subset\Lambda$ acts. We will show that
$H^\bullet_{T\times\Gm}(\Gr_G)\cong\kk[\Lambda]\otimes_{\kk[Q]}
\O\left(\on{Bl}^{\on{aff}}_\Gamma(\t_{\on{aff}}\times_{\BA^1}(\t_{\on{aff}}/W))\right)\cong\bigoplus_{\pi_1(G)}
\O\left(\on{Bl}^{\on{aff}}_\Gamma(\t_{\on{aff}}\times_{\BA^1}(\t_{\on{aff}}/W))\right)$
(note that $\pi_1(G)=\Lambda/Q$). It will follow that
$H^*_{\GO \rtimes \Gm}(\Gr_G)\cong\Big(\kk[\Lambda]\otimes_{\kk[Q]}\O\left(\on{Bl}^{\on{aff}}_\Gamma(\t_{\on{aff}}\times_{\BA^1}(\t_{\on{aff}}/W))\right)\Big)^W\cong\bigoplus_{\pi_1(G)}\O(N_{\Lt^*/W\times \Lt^*/W}\Delta)$.

Let $G^{\mathrm{sc}}$ stand for the simply connected cover of $G$. Then the
Grassmannian $\Gr_G$ is a union of connected components indexed by 
$\chi\in\Lambda/Q$, and each connected component is isomorphic
to $\Gr_{G^{\mathrm{sc}}}$. More precisely, we have a natural action of the coweight lattice
$\Lambda=X_*(T)\subset{\mathbf G}_{\mathbf F}$ on $\Gr_G$. Let us pick a coweight
$\lambda\in\Lambda$ sending the neutral component $\Gr_{G^{\mathrm{sc}}}$ to another component
$\Gr_G^\chi$. Then we can $W$-equvariantly identify $H^\bullet_{T\times\Gm}(\Gr_{G^{\mathrm{sc}}})$
with $H^\bullet_{T\times\Gm}(\Gr_G^\chi)$ in the following way. The desired identification
is the composition of the isomorphism induced by (the shift by) $\lambda$ on $\Gr_G$,
with the inverse of the action of $\lambda$ on $\on{Bl}^{\on{aff}}_\Gamma(\t_{\on{aff}}\times_{\BA^1}(\t_{\on{aff}}/W))$.
It is easy to see that the resulting isomorphism is independent of the choice of $\lambda$ and hence it is
$W$-equivariant.

Now $H^\bullet_{T\times\Gm}(\Gr_G)$ is equal to $\bigoplus_{\chi\in\Lambda/Q}H^\bullet_{T\times\Gm}(\Gr_G^\chi)$. Finally,
$H^\bullet_{\GO\rtimes\Gm}(\Gr_G)=H^\bullet_{G\times\Gm}(\Gr_G)$ coincides
with the $W$-invariants in $H^\bullet_{T\times\Gm}(\Gr_G)$.
\hfill $\Box$

\subsection{Canonical filtration on $H^\bullet_{T\times\Gm}(\Gr_G,\CF)$}
\label{can}
For a $\GO\rtimes\Gm$-equivariant perverse sheaf $\CF$ on $\Gr_G$ we will
define a canonical filtration on $H^\bullet_{T\times\Gm}(\CF)=
H^\bullet_{\GO\rtimes\Gm}(\CF)\otimes_{\O(\t/W)}\O(\t)$. Note that
if $\pi$ is the projection $\t\to\t/W$, and $(\pi,\Id,\Id)$ is the 
projection $\t\times(\t/W)\times\Aone\to\t/W\times\t/W\times\Aone$,
then $H^\bullet_{\GO\rtimes\Gm}(\CF)\otimes_{\O(\t/W)}\O(\t)=
(\pi,\Id,\Id)^*H^\bullet_{\GO\rtimes\Gm}(\CF)$.
For $\lambda\in X_*(T)$ we denote by $\lambda$ the corresponding $T$-fixed
point of $\Gr_G$. We denote by $\fT_\lambda$ the semiinfinite 
$N_-({\mathbf F})$-orbit through $\lambda$.
We denote by $\ofT_\lambda$ the closure of $\fT_\lambda$, that is the union
of $\fT_\mu$ over $\mu\geq\lambda$. We filter
$H^\bullet_{T\times\Gm}(\CF)$ by the images of
$r_\lambda:\ H^\bullet_{\ofT_\lambda,T\times\Gm}(\Gr_G,\CF)\to
H^\bullet_{T\times\Gm}(\Gr_G,\CF)$ (cohomology with supports).
The associated graded of this filtration is 
$\bigoplus_\lambda H^\bullet_{\fT_\lambda,T\times\Gm}(\Gr_n\CF)$.
Since $\lambda$ is the only $T\times\Gm$-fixed point of $\fT_\lambda$, we have
$H^\bullet_{\fT_\lambda,T\times\Gm}(\Gr_G,\CF)=H^\bullet_{T\times\Gm}(\lambda)
\otimes\fj_\lambda^*\imath_\lambda^!\CF$ where $\imath_\lambda$ is the
locally closed embedding of $\fT_\lambda$ into $\Gr_G$, and 
$\fj_\lambda$ is the embedding of $\lambda$ into $\fT_\lambda$.

Now recall that $\CF\mapsto\fj_\lambda^*\imath_\lambda^!\CF$ is the
$\lambda$-weight component of the Mirkovi\'c-Vilonen fiber functor on
the category of $\GO\rtimes\Gm$-equivariant perverse sheaves on $\Gr_G$.
In other words, if $V$ is a $\LG$-module, and $\CF=S(V)$, then
$\CF\mapsto\fj_\lambda^*\imath_\lambda^!\CF=\ _\lambda V$ where
$_\lambda V$ is the $\lambda$-weight component of $V$.

Furthermore, we claim that the $\O(\t\times(\t/W)\times\Aone)$-module
$H^\bullet_{T\times\Gm}(\lambda)$ is canonically isomorphic to
$(\Id,\pi,Id)_*\O(\Gamma_\lambda)$ where 
$\Gamma_\lambda\subset\t\times\t\times\Aone$ is given by the equation
$\Gamma_\lambda=\{(x_1,x_2,a):\ x_2=x_1+a\lambda\}$. In effect, 
let $p$ stand for the projection from the affine flag variety $\Fl_G$ to
the affine Grassmannian $\Gr_G$. Let $\tilde\lambda$ be a 
$T\times\Gm$-fixed point of $\Fl_G$ such that $p$ projects $\tilde\lambda$
isomorphically onto $\lambda$. Let $\bI$ stand for the Iwahori subgroup of
$\GO$. We have $H^\bullet_{T\times\Gm}(\Fl_G)=
H^\bullet_\Gm(\bI\backslash\GF/\bI)$, and so 
$H^\bullet_{T\times\Gm}({\tilde\lambda})$
is a module over $\O(\t\times\t\times\Aone)$. Clearly, the 
$\O(\t\times(\t/W)\times\Aone)$-module $H^\bullet_{T\times\Gm}(\lambda)$ is
isomorphic to the direct image of the 
$\O(\t\times\t\times\Aone)$-module $H^\bullet_{T\times\Gm}({\tilde\lambda})$
under the projection $\t\times\t\times\Aone\to\t\times(\t/W)\times\Aone$.
So it suffices to check that the
$\O(\t\times\t\times\Aone)$-module $H^\bullet_{T\times\Gm}({\tilde\lambda})$
is isomorphic to $\O(\Gamma_\lambda)$ after localization along $\t\times\Aone$.

The set of $T$-fixed points in $\Fl_G$ is canonically identified with the
extended affine Weyl group $W_{aff}$ of $G$, and we choose $\tilde\lambda$
so that it coincides 
with $\lambda\in W_{aff}$. Then the preimage $T_\lambda$ of 
${\tilde\lambda}\in\Fl_G$ in $\GF$ is homotopically equivalent to 
$T$, and the action of $T\times T\times\Gm$ on $T_\lambda$ is homotopically
equivalent to $(t_1,t_2,z)(t)=t_1tt_2^{-1}\lambda(z)$. We conclude that the
$\O(\t\times\t\times\Aone)$-module $H^\bullet_{T\times\Gm}({\tilde\lambda})=
H^\bullet_\Gm(T\backslash T_\lambda/T)$
is isomorphic to $\O(\Gamma_\lambda)$.

We have proved the following

\begin{lem}
\label{Harvard}
For $V\in Rep(\LG)$, the $\O(\t\times(\t/W)\times\Aone)$-module
$H^\bullet_{T\times\Gm}(\Gr_G,S(V))$ has a canonical filtration with
the associated graded $\bigoplus_\lambda (\Id,\pi,\Id)_*
\O(\Gamma_\lambda)\otimes\ _\lambda V$. In particular, 
$H^\bullet_{T\times\Gm}(\Gr_G,S(V))$ is flat as an $\O(\t\times\Aone)$-module.
\end{lem}
\hfill $\Box$

\subsection{Levi-equivariant cohomology}
\label{levi}
Let $T\subset L\subset G$ be a Levi subgroup. We denote by $P_L$ (resp.
$P^-_L$) the parabolic subgroup generated by $L$ and the positive (resp.
negative) Borel subgroup $B$ (resp. $B_-$). We denote by $W_L\subset W$ 
the Weyl group of $L$. We denote by $\pi_L$ the projection from 
$\t/W_L$ to $\t/W$. We denote by $X^+_L$ the set of highest weights of
irreducible $\LL$-modules, where $\LL\subset\LG$ stands for the Langlands
dual Levi subgroup. We have a natural projection from $X^+_L$ to the lattice
$X^*(Z(\LL))$ of characters of the center $Z(\LL)$ of $\LL$.
The set of $[P^-_L,P^-_L]({\mathbf F})$-orbits in $\Gr_G$ is indexed by
$X^*(Z(\LL))$. For $\lambda\in X^*(Z(\LL))$ we will denote the corresponding
orbit by $_L\fT_\lambda$, and its closure by $_L\ofT_\lambda$.
The locally closed embedding of $\fT_\lambda$ into $\Gr_G$ is 
denoted by $\imath_\lambda$.
For an $\LL$-module $V$ we denote by $S_L(V)$ the corresponding
$L({\mathbf O})\rtimes\Gm$-equivariant perverse sheaf on $\Gr_L$.

\begin{lem}
\label{Chicago}
For $V\in Rep(\LG)$, the $\O(\t/W_L\times\t/W\times\Aone)$-module
$(\pi_L,\Id,\Id)^*H^\bullet_{\GO\rtimes\Gm}(\Gr_G,S(V))$
carries a canonical filtration $F^\bullet_L$ such that the associated graded is
equipped with a canonical isomorphism
$$\ ^{\operatorname{top}}\Xi_L:\ gr(\pi_L,\Id,\Id)^*H^\bullet_{\GO\rtimes\Gm}(\Gr_G,S(V))
\iso(\Id,\pi_L,\Id)_*
H^\bullet_{L({\mathbf O})\rtimes\Gm}(\Gr_L,S_L(V|_{\LL})).$$
\end{lem}

\proof: We have 
$(\pi_L,\Id,\Id)^*H^\bullet_{\GO\rtimes\Gm}(\Gr_G,S(V))=
H^\bullet_{L({\mathbf O})\rtimes\Gm}(\Gr_G,S(V))=H^\bullet_{L\times\Gm}
(\Gr_G,S(V))$.
The canonical filtration in question is filtration by the images of
$r_\lambda:\ H^\bullet_{_L\ofT_\lambda,L\times\Gm}(\Gr_G,S(V))\to
H^\bullet_{L\times\Gm}(\Gr_G,S(V))$ (cohomology with supports; here
$\lambda\in X^*(Z(\LL))$). The associated graded of this filtration is
$\bigoplus_{\lambda\in X^*(Z(\LL))}H^\bullet_{_L\fT_\lambda,L\times\Gm}
(\Gr_G,S(V))$. Let $\fp^\lambda:\ _L\fT_\lambda\to\Gr_L$ denote the
natural $L({\mathbf O})\rtimes\Gm$-equivariant projection. Then we have
$H^\bullet_{_L\fT_\lambda,L\times\Gm}(\Gr_G,S(V))=
H^\bullet_{L\times\Gm}(\Gr_L,\fp^\lambda_*\imath_\lambda^!S(V))$.
However, according to~\cite{BD}, we have a canonical isomorphism
$\bigoplus_{\lambda\in X^*(Z(\LL))}\fp^\lambda_*\imath_\lambda^!S(V)=
S_L(V|_{\LL})$. The lemma is proved.
\hfill $\Box$

\subsection{Transitivity for a pair of Levi subgroups}
\label{toppair}
We have a canonical isomorphism
$^{\operatorname{top}}\Xi_L:\ 
gr(\pi_L,\Id,\Id)^*H^\bullet_{\GO\rtimes\Gm}(\Gr_G,S(V))
\iso(\Id,\pi_L,\Id)_*
H^\bullet_{L({\mathbf O})\rtimes\Gm}(\Gr_L,S_L(V|_{\LL}))$. 
In the RHS we have the restriction of 
$\O(\t/W_L\times\t/W_L\times\Aone)$-module
$H^\bullet_{L({\mathbf O})\rtimes\Gm}(\Gr_L,S_L(V|_{\LL}))$ to
$\O(\t/W_L\times\t/W\times\Aone)$. To save a bit of notation in what 
follows we will write simply
$$^{\operatorname{top}}\Xi_L:\ 
gr(\pi_L,\Id,\Id)^*H^\bullet_{\GO\rtimes\Gm}(\Gr_G,S(V))
\iso H^\bullet_{L({\mathbf O})\rtimes\Gm}(\Gr_L,S_L(V|_{\LL})).$$

If $\LT\subset\LL'\subset\LL$
is another Levi subgroup, then we denote by $\pi_{L'}^L$ the projection
from $\t/W_{L'}$ to $\t/W_L$. Note that the filtration $F_{L'}^\bullet$
on $(\pi_{L'},\Id,\Id)^*H^\bullet_{\GO\rtimes\Gm}(\Gr_G,S(V))=
(\pi_{L'}^L,\Id,\Id)^*(\pi_L,\Id,\Id)^*H^\bullet_{\GO\rtimes\Gm}(\Gr_G,S(V))$
is a refinement of the filtration $(\pi_{L'}^L,\Id,\Id)^*F_L^\bullet$,
and hence induces a canonical filtration $F_{L'}^{L\bullet}$ on
$(\pi_{L'}^L,\Id,\Id)^*gr_{F_L^\bullet}(\pi_L,\Id,\Id)^*
H^\bullet_{\GO\rtimes\Gm}(\Gr_G,S(V))$. The isomorphism
$(\pi_{L'}^L,\Id,\Id)^*\ ^{\operatorname{top}}\Xi_L$ carries the filtration $F_{L'}^{L\bullet}$
to the filtration $F_{L'}^\bullet$ on
$(\pi_{L'}^L,\Id,\Id)^*
H^\bullet_{L({\mathbf O})\rtimes\Gm}(\Gr_L,S_L(V|_{\LL}))$.
We have a canonical isomorphism
$\ ^{\operatorname{top}}\Xi_{L'}^L:\ gr_{F_{L'}^\bullet}(\pi_{L'}^L,\Id,\Id)^*
H^\bullet_{L({\mathbf O})\rtimes\Gm}(\Gr_L,S_L(V|_{\LL}))\iso
H^\bullet_{L'({\mathbf O})\rtimes\Gm}(\Gr_{L'},S_{L'}(V|_{\LL'}))$.
We consider the composition
\begin{equation}
\begin{CD} 
gr_{F_{L'}^\bullet}(\pi_{L'},\Id,\Id)^*H^\bullet_{\GO\rtimes\Gm}(\Gr_G,S(V))
@>gr_{F_{L'}^{L\bullet}}(\pi_{L'}^L,\Id,\Id)^*\ ^{\operatorname{top}}\Xi_L>>
\end{CD}
\end{equation}
$$\begin{CD}
gr_{F_{L'}^\bullet}(\pi_{L'}^L,\Id,\Id)^*
H^\bullet_{L({\mathbf O})\rtimes\Gm}(\Gr_L,S_L(V|_{\LL}))
@>\ ^{\operatorname{top}}\Xi_{L'}^L>>
H^\bullet_{L'({\mathbf O})\rtimes\Gm}(\Gr_{L'},S_{L'}(V|_{\LL'}))
\end{CD}$$
Then we have
\begin{lem}
\label{tp}
$^{\operatorname{top}}\Xi_{L'}^L\circ gr_{F_{L'}^{L\bullet}}
(\pi_{L'}^L,\Id,\Id)^*\ ^{\operatorname{top}}\Xi_L=\ 
^{\operatorname{top}}\Xi_{L'}$.
\end{lem}
\hfill $\Box$

\subsection{Tensor structure on equivariant cohomology}
\label{tenscoh}
A $\GO\rtimes\Gm$-equivariant sheaf $\CF$ on $\Gr_G$ will be viewed as
a sheaf on the stack $\GO\rtimes\Gm\backslash\GF\rtimes\Gm/\GO\rtimes\Gm$.
Given two such sheaves $\CF_1,\CF_2$ we define 
$\CF_1\boxtimes_{\GO\rtimes\Gm}\CF_2$ as the descent of
$\CF_1\boxtimes\CF_2$ from 
$$\left(\GO\rtimes\Gm\backslash\GF\rtimes\Gm/\GO\rtimes\Gm\right)\times
\left(\GO\rtimes\Gm\backslash\GF\rtimes\Gm/\GO\rtimes\Gm\right)$$ to
$\GO\rtimes\Gm\backslash(\GF\rtimes\Gm\times_{\GO\rtimes\Gm}\GF\rtimes\Gm)/
\GO\rtimes\Gm$. Clearly, $$H^\bullet(
\GO\rtimes\Gm\backslash(\GF\rtimes\Gm\times_{\GO\rtimes\Gm}\GF\rtimes\Gm)/
\GO\rtimes\Gm,\CF_1\boxtimes_{\GO\rtimes\Gm}\CF_2)=$$
$$H^\bullet(\GO\rtimes\Gm\backslash\GF\rtimes\Gm/\GO\rtimes\Gm,\CF_1)
\otimes_{H^\bullet_{\GO\rtimes\Gm}(pt)}
H^\bullet(\GO\rtimes\Gm\backslash\GF\rtimes\Gm/\GO\rtimes\Gm,\CF_2)=$$
$$H^\bullet_{\GO\rtimes\Gm}(\Gr_G,\CF_1)\otimes_{\O(\t/W\times\Aone)}
H^\bullet_{\GO\rtimes\Gm}(\Gr_G,\CF_2)=:
H^\bullet_{\GO\rtimes\Gm}(\Gr_G,\CF_1)\star
H^\bullet_{\GO\rtimes\Gm}(\Gr_G,\CF_2).$$
The multiplication in $\GF\rtimes\Gm$ gives rise to the map
$$m:\ \GO\rtimes\Gm\backslash(\GF\rtimes\Gm\times_{\GO\rtimes\Gm}
\GF\rtimes\Gm)/
\GO\rtimes\Gm\to\GO\rtimes\Gm\backslash\GF\rtimes\Gm/\GO\rtimes\Gm.$$
The convolution $\CF_1*\CF_2$ is defined as
$m_*(\CF_1\boxtimes_{\GO\rtimes\Gm}\CF_2)$. Hence
$$H^\bullet_{\GO\rtimes\Gm}(\Gr_G,\CF_1*\CF_2)=
H^\bullet(\GO\rtimes\Gm\backslash\GF\rtimes\Gm/\GO\rtimes\Gm,\CF_1*\CF_2)=$$
$$H^\bullet(
\GO\rtimes\Gm\backslash(\GF\rtimes\Gm\times_{\GO\rtimes\Gm}\GF\rtimes\Gm)/
\GO\rtimes\Gm,\CF_1\boxtimes_{\GO\rtimes\Gm}\CF_2)=$$
$$H^\bullet_{\GO\rtimes\Gm}(\Gr_G,\CF_1)\star
H^\bullet_{\GO\rtimes\Gm}(\Gr_G,\CF_2).$$
Now for $V_1,V_2\in Rep(\LG)$, and $\CF_1=S(V_1),\ \CF_2=S(V_2)$ we have
a canonical isomorphism $S(V_1\otimes V_2)\iso S(V_1)*S(V_2)$, and thus
$$\ ^{\operatorname{top}}\omega_{V_1,V_2}:\ 
H^\bullet_{\GO\rtimes\Gm}(\Gr_G,S(V_1\otimes V_2))\iso
H^\bullet_{\GO\rtimes\Gm}(\Gr_G,S(V_1))\star
H^\bullet_{\GO\rtimes\Gm}(\Gr_G,S(V_2)).$$
According to Lemma~\ref{Harvard} (cf. also Lemma~\ref{Chicago}),
we have a canonical isomorphism 
$$^{\operatorname{top}}\Xi_V=\ ^{\operatorname{top}}\Xi_{T,V}:\ 
gr(\pi,\Id,\Id)^*H^\bullet_{\GO\rtimes\Gm}(\Gr_G,S(V))\iso
(\Id,\pi,\Id)_*H^\bullet_{T\times\Gm}(\Gr_T,S_T(V|_{\LT})).$$ 
In the RHS we have the restriction of
$\O(\t\times\t\times\Aone)$-module $H^\bullet_{T\times\Gm}(\Gr_T,
S_T(V|_{\LT}))$ to
$\O(\t\times(\t/W)\times\Aone)$. To save a bit of notation in what
follows we will write simply 
$\ ^{\operatorname{top}}\Xi_V:\ 
gr(\pi,\Id,\Id)^*H^\bullet_{\GO\rtimes\Gm}(\Gr_G,S(V))\iso
H^\bullet_{T\times\Gm}(\Gr_T,S_T(V|_{\LT}))$. 
It follows that after tensoring with
$\kk(\t\times\Aone)$ (over the first and third factors in
$\O(\t/W\times\t/W\times\Aone)$) we have a canonical isomorphism
$$\ ^{\operatorname{top}}\Xi_V^\gen:\ H^\bullet_{\GO\rtimes\Gm}(\Gr_G,S(V))
\otimes_{\O(\t/W\times\Aone)}\kk(\t\times\Aone)=$$
$$=gr(\pi,\Id,\Id)^*H^\bullet_{\GO\rtimes\Gm}(\Gr_G,S(V))
\otimes_{\O(\t\times\Aone)}\kk(\t\times\Aone)\iso$$
$$\iso H^\bullet_{T\times\Gm}(\Gr_T,S_T(V|_{\LT}))
\otimes_{\O(\t\times\Aone)}\kk(\t\times\Aone)=
\bigoplus_\lambda\left(\O(\Gamma_\lambda)
\otimes_{\O(\t\times\Aone)}\kk(\t\times\Aone)\right)\otimes\ 
_\lambda V$$
Now we have a canonical isomorphism 
$\O(\Gamma_\mu)\star\O(\Gamma_\nu):=\O(\Gamma_\mu)
\otimes_{\O(\t\times\Aone)}\O(\Gamma_\nu)=\O(\Gamma_{\mu+\nu})$.
Hence we get a canonical isomorphism 
$$\ ^{\operatorname{top}}\Xi_{V_1}^\gen\star\ 
^{\operatorname{top}}\Xi_{V_2}^\gen:\
(H^\bullet_{\GO\rtimes\Gm}(\Gr_G,S(V_1))\star
H^\bullet_{\GO\rtimes\Gm}(\Gr_G,S(V_2))) 
\otimes_{\O(\t/W\times\Aone)}\kk(\t\times\Aone)\iso$$
$$\iso\bigoplus_{\mu+\nu=\lambda}\left(\O(\Gamma_\lambda)
\otimes_{\O(\t\times\Aone)}\kk(\t\times\Aone)\right)\otimes\
_\mu V_1\otimes\ _\nu V_2=$$
$$=\bigoplus_\lambda\left(\O(\Gamma_\lambda)
\otimes_{\O(\t\times\Aone)}\kk(\t\times\Aone)\right)\otimes\
_\lambda(V_1\otimes V_2)$$
We want to compare it with $\ ^{\operatorname{top}}\Xi_{V_1\otimes V_2}^\gen:$
$$H^\bullet_{\GO\rtimes\Gm}(\Gr_G,
S(V_1\otimes V_2))\otimes_{\O(\t/W\times\Aone)}\kk(\t\times\Aone)
\iso\bigoplus_\lambda\left(\O(\Gamma_\lambda)
\otimes_{\O(\t\times\Aone)}\kk(\t\times\Aone)\right)\otimes\
_\lambda(V_1\otimes V_2)$$

\begin{prop}
\label{tenstop}
$^{\operatorname{top}}\Xi_{V_1\otimes V_2}^\gen=
(\ ^{\operatorname{top}}\Xi_{V_1}^\gen\star\ 
^{\operatorname{top}}\Xi_{V_2}^\gen)
\circ\ ^{\operatorname{top}}\omega_{V_1,V_2}$.
\end{prop}

\proof: The equality readily reduces to the following compatibility.
Let $\Phi_{MV}=\oplusl_\la \Phi_{MV}^\lambda :\CF\mapsto \oplusl_\lambda 
\fj_\lambda^*\imath_\lambda^!\CF$ be the
 Mirkovi\'c-Vilonen fiber functor on the Satake category
$Perv_{\GO}(\Gr)$ (notations of~\ref{can}, see~\cite{MV}). 
The (proof of) Lemma 1 provides
a canonical isomorphism 
\begin{equation}\label{isoten}
gr(H^\bu_{T\times \Gm}(Gr_G,\F) )\cong \oplusl 
\Phi_{MV}^\la\otimes \CO(\Gamma_\la).
\end{equation}
We have to check that this isomorphism is compatible with the tensor
structure, i.e. for $\CF,\CG\in Perv_{\GO}(\Gr)$ we have to check coincidence
of the two embeddings from $\Phi_{MV}(\CF)\otimes \Phi_{MV}(\CG)$
in $gr(H^\bu_{T\times \Gm}(Gr_G,\F*\CG) )$, where the first one
comes from the isomorphisms \eqref{isoten} for $\CF$, $\CG$ and tensor
structure on the functor $gr(H^\bu_{T\times \Gm})$, and the second one
comes from the tensor structure on the functor $\Phi_{MV}$ and
isomorphism \eqref{isoten} for $\CF*\CG$.

To check the equality, we recall a ``filtration'' in the $I$-equivariant
 derived category
on a $\GO$ equivariant perverse sheaf $\F$, which induces the above filtration
on $H^\bu_{T\times \Gm}(Gr_G,\F)$. Here and below by a ``filtration''
on an object $X$ of a triangulated category we mean a collection
of object $X_0=0,\dots, X_n=X$ and distinguished triangles
$X_i\to X_{i+1} \to Y_i$; the objects $Y_i$ will be called the ``subquotients''
of the filtration.

 Let $\CalC$, $\CalD$ be  the equivariant constructible derived category
with respect to the natural action of $\bI\rtimes\Gm$ on $\Gr_G$ and
 on $\Fl_G$ respectively.
Thus convolution $*_I$ provides
$\CalD$ with a monoidal structure, and $\CalC$
with an action of the monoidal  category $\CalD$.

Recall the {\em Wakimoto sheaves} $J_\la\in \CalD$, characterized by the
following properties: $J_{\la+\mu}\cong J_\la *_I J_\mu$, while for
 a dominant weight $\la \in \Lambda^+$ the sheaf
$J_\la$ is the $*$-extension
of the constant perverse sheaf from the Iwahori orbit corresponding to $\la$,
 see, e.g. \cite{AB}.

%{\bf though there may be a big normalization issue here}

Recall that $p$ stands for the projection $\Fl_G\to\Gr_G$.
We set $J_\la^{\Gr}=p_*(J_\la)$. It is not hard to show
that for $\F\in D_{\GO\ltimes \Gm}(\Gr)$, $J_\la^{\Gr}*\F\cong J_\la*_I \F$
canonically. Also $J^{\Gr}_\la$ can be characterized by 
$\fj_\mu^*\imath_\mu^!J_\la^{\Gr} = \kk^{\delta_{\la,\mu}}$
(cf. \cite{AB}).

Fix $\F\in Perv_ {\GO\ltimes \Gm}(\Gr)$, and choose a coweight 
$\la$ deep inside the
dominant chamber. Then one shows that $j_\nu^!(J_\la*\F)\cong 
\kk[\dim Gr_\nu] \otimes \fj_{\nu-\la}^* \imath_{\nu-\la}^!\F$
for all $\nu$. Thus one can consider the Cousin ``filtration'' 
on $J_\la*\F$ with subquotients $j_{\nu*}j_\nu^! (J_\la*\F)$ and
apply the functor $J_{-\la}*_I\ \ $ to it, thereby obtaining
a ``filtration'' on $\CF$  with ``subquotients"
$\Phi_{MV}^\mu(\F)\otimes J_\mu^{\Gr}$. It is clear that this
``filtration'' induces the above filtration on $H_{T\times \Gm}(\F)$.

Let now $\CF,\CG$ be a pair of objects of $Perv_{\GO\ltimes \Gm}(\Gr)$.
The above ``filtration'' on $\CF$ induces a "filtration" on $\CF*\CG$
with ``subquotients'' 
$$\Phi_{MV}^\mu(\F)\otimes J_\mu^{\Gr} *\CG=
\Phi_{MV}^\mu(\F)\otimes J_\mu*_I \CG.$$
Using the ``filtration'' on $\CG$ with ``subquotients'' 
$\Phi_{MV}^\nu(\CG)\otimes
J_\nu^{\Gr}$ we get a ``filtration'' on $\CF*\CG$ with ``subquotients''
$$\Phi_{MV}^\mu(\F)\otimes \Phi_{MV}^\nu(\G) \otimes J_\mu*_I J_\nu^{\Gr}=
\Phi_{MV}^\mu(\F)\otimes \Phi_{MV}^\nu(\G) \otimes J_{\mu+\nu}^{\Gr}.$$
Comparing it with the ``filtration'' on $\CF*\CG$
with ``subquotients'' $\Phi_{MV}^\eta (\F*\G)$ we get an isomorphism
$\Phi_{MV}(\CF*\CG)\cong \Phi_{MV}(\CF)\otimes \Phi_{MV}(\CG)$.
It is not hard to see that this isomorphism coincides with any of the standard
definitions of tensor structure on $\Phi_{MV}$; in fact, a close
description of the tensor structure appears in \cite{BG},  Theorem 3.2.8.

Now we see that the isomorphism 
$$gr\, H_{T\times \Gm}^\bu(\CF*\CG ) \cong 
gr\, H_{T\times \Gm}^\bu(\CF)\star\ 
gr\, H_{T\times \Gm}^\bu(\CG )$$
breaks as a direct sum of maps 
$$\left(\Phi_{MV}^\mu(\CF) \otimes 
H_{T\times \Gm}^\bu(J_\mu^{\Gr})\right)\star\ 
\left(\Phi_{MV}^\nu(\CG) \otimes H_{T\times \Gm}^\bu(J_\nu^{\Gr})\right)
\to \Phi_{MV}^{\mu+\nu}
(\CF*\CG) \otimes  H_{T\times \Gm}^\bu(J_{\mu+\nu}^{\Gr}),
$$
coming from the map $\Phi_{MV}^\mu(\CF)\otimes \Phi_{MV}^\nu(\CG)
\to \Phi_{MV}^{\mu+\nu}(\CF*\CG)$ induced by the tensor structure
on $\Phi_{MV}$, and the natural isomorphism
$$ \CO(\Gamma_\mu)\star\CO(\Gamma_\nu)=
H_{T\times \Gm}^\bu(J_\mu^{\Gr})\star
 H_{T\times \Gm}^\bu(J_\nu^{\Gr})
 \cong  H_{T\times \Gm}^\bu(J_{\mu+\nu}^{\Gr})=\CO(\Gamma_{\mu+\nu}).$$
The claim follows. \hfill $\Box$

\section{Algebra}

\subsection{Some properties of the Kostant functor $\kappa_\hbar$} 
The following properties of the Kostant functor will play an important role
in the proof of the main results.

\begin{lem}
\label{propKons}
a) The functors $\kappa,\, \kappah$ are exact.

b) The functors %$\kappa|_{Coh_{fr}^\LG(\Lg^*)}$,  
$\kappa|_{Coh_{fr}^{\LG\times \Gm}(\Lg^*)}$,
$\kappah|_{\HChfr}$ are full embeddings.
%\footnote{This is true for all groups, not necessarily
%adjoint, however strange this may seem.}
\end{lem} 

\proof is essentially due to B.~Kostant, cf.~\cite{K1}.

a) For the exactness of $\kappa$, note that the functor of restriction
$\CF\mapsto\CF|_\Upsilon$ is exact on $Coh^\LG(\Lg^*)$, and then
$\CF|_\Upsilon$ is an $\LN_-$-equivariant coherent sheaf on $\Upsilon$.
Recall that $\LN_-$ acts on $\Upsilon$ freely, and
$\Upsilon/\LN_-\simeq\Sigma$. Hence the functor of invariants
$\CG\mapsto\CG^{\LN_-}$ is exact on $Coh^{\LN_-}(\Upsilon)$.
The exactness of $\kappa$ follows.

For the exactness of $\kappa_\hbar$, we will prove that both the 
functors of $-\psi$-coinvariants, and $\LN_-$-invariants are exact,
and hence $\kappa_\hbar$ is exact as their composition.
It is enough to check it on the 
positively graded $\hbar$-Harish-Chandra bimodules. Then it is enough
to check the exactness on the subcategory of $\hbar$-Harish-Chandra
bimodules with grading degrees between 0 and $n$ for a fixed $n\gg0$.
Thus it suffices to check the exactness on the subcategory of 
$\hbar$-Harish-Chandra bimodules with nilpotent action of $\hbar$,
and then it suffices to consider the subcategory of
$\hbar$-Harish-Chandra bimodules with the {\em trivial} action of
$\hbar$. However, an $\hbar$-Harish-Chandra bimodule $M$ with the
trivial action of $\hbar$ is nothing else than a 
$\LG$-equivariant coherent sheaf
on $\Lg^*$, and $M\Lotimes_{U_\hbar(\Ln_-)_2}(-\psi)=M|_\Upsilon$.
In particular, the functor of $-\psi$-coinvariants is exact according
to the previous paragraph. For the same reason,
$M\mapsto(M\Lotimes_{U_\hbar(\Ln_-)_2}(-\psi))^{\LN_-}=(M|_\Upsilon)^{\LN_-}$
is exact. This completes the proof of a).

b) $\kappa|_{Coh_{fr}^{\LG\times \Gm}(\Lg^*)}$ is fully faithful
since the complement to $(\Lg^*)^{reg}$ in $\Lg^*$ has codimension~2,
and the centralizer of a generic regular element is connected.

To prove that $\kappah|_{\HChfr}$ is fully faithful, we consider
free $\hbar$-Harish-Chandra bimodules $M_1=U_\hbar\otimes V_1,\
M_2=U_\hbar\otimes V_2$, and the following commutative diagram:

\begin{equation}
\begin{CD}
\Hom(M_1,M_2)@>\epsilon>>\Hom(\kappah M_1,\kappah M_2)@>\delta>>
\Hom(\kappah M_1,\kappah M_2)/\hbar\\
@VV\beta V @. @VV\gamma V\\
\Hom(M_1/\hbar,M_2/\hbar)@>\alpha>>\Hom(\kappa(M_1/\hbar),\kappa(M_2/\hbar))
@>\sim>>\Hom((\kappah M_1)/\hbar,(\kappah M_2)/\hbar)
\end{CD}
\end{equation}

We have just proved that $\alpha$ is an isomorphism. Moreover, $\beta$
is surjective since $\Hom(U_\hbar\otimes V_1,U_\hbar\otimes V_2)=
\Hom_\LG(V_1\otimes V_2^*,U_\hbar)$, and all the $\LG$-modules in
question are semisimple. It follows that $\gamma$ is surjective.
On the other hand, $\gamma$ is injective since $\kappah M_1,\ \kappah
M_2$ are free over $\kk[\hbar]$. Now that $\gamma$ is proved to be an
isomorphism, the composition $\delta\circ\epsilon$ must be surjective.
Hence $\epsilon$ is surjective by Nakayama Lemma. It remains to prove
that $\epsilon$ is injective. Since $\kappah$ is exact, it is enough
to prove that $\kappah M\ne0$ for a nonzero subobject $M$ of a free
$\hbar$-Harish-Chandra bimodule $M_2$. We consider a nonzero subobject
$M/\hbar\subset M_2/\hbar$ of a free $\CO(\Lg^*)$-module $M_2/\hbar$.
It suffices to prove that $\kappa M\ne0$. However, the support of any
nonzero section of a free $\CO(\Lg^*)$-module is the whole of $\Lg^*$,
hence its restriction to $\Upsilon$ is nonzero.

The lemma is proved.
\hfill $\Box$

\subsection{De-symmetrized Kostant functor $\varkappa_\hbar$} 
We denote by $\pi$ the projection $\Lt^*\to\Lt^*/W$, and
we denote by $(\pi,\Id,\Id)$ the projection $\Lt^*\times(\Lt^*/W)\times\Aone
\to\Lt^*/W\times\Lt^*/W\times\Aone$. For $V\in Rep(\LG)$ we are going
to describe $(\pi,\Id,\Id)^*\phi(V)\in Coh(\Lt^*\times(\Lt^*/W)\times\Aone)$.

To this end we consider the universal Verma module 
$\CM_\hbar(-\rho)=
U_\hbar\otimes_{U_\hbar(\Lb)}\kk[\hbar][\Lt](-\rho)$,
and $\kk[\hbar][\Lt](-\rho)$ is a $U_\hbar(\Lb)$-module which factors
through the $U_\hbar(\Lt)=\kk[\hbar][\Lt]$-module where $t\in\Lt$ acts by
multiplication by $t-\hbar\rho(t)$ (recall that $\rho$ is the halfsum
of positive roots of $\Lg$). For an
$\hbar$-Harish-Chandra bimodule $M$ we set
$\varkappa_\hbar(M):=\CM_\hbar(-\rho)\overset{L}{\otimes}_{U_\hbar}M
\overset{L}{\otimes}_{U_\hbar(\Ln_-)_2}\psi=
\kk[\hbar][\Lt](-\rho)\overset{L}{\otimes}_{U_\hbar(\Lb)_1}M
\overset{L}{\otimes}_{U_\hbar(\Ln_-)_2}\psi$. This is an
$\O(\Lt^*\times(\Lt^*/W)\times\Aone)$-module: the action of 
$\O(\Lt^*\times\Aone)=U_\hbar(\Lt)$ comes from the fact that $U_\hbar(\Lt)$
normalizes $U_\hbar(\Lb)$, and the action of $\O(\Lt^*/W\times\Aone)$ is the
action of the center $Z(U_\hbar)$ of the second copy of $U_\hbar$, as before.

For $V\in Rep(\LG)$ we set $\varphi(V):=\varkappa_\hbar(Fr(V))$.

\begin{lem}
\label{mit}
For $V\in Rep(\LG)$ we have a canonical isomorphism
$\varphi(V)\simeq(\pi,\Id,\Id)^*\phi(V)$.
\end{lem}

\proof: We denote by $\CW_\hbar^-:=U_\hbar\otimes_{U_\hbar(\Ln_-)}\psi$ 
the Whittaker
$U_\hbar$-module. By a theorem of Kostant, 
$\End_{U_\hbar}(\CW_\hbar^-)=Z(U_\hbar)$, and the category $\CA$ of 
$(U_\hbar(\Ln_-),\psi)$-integrable $U_\hbar$-modules is equivalent to the
category of $Z(U_\hbar)$-modules (here a $U_\hbar(\Ln_-)$-module $L$
is called $(U_\hbar(\Ln_-),\psi)$-integrable if the action of $\Ln_-$
on $L\otimes(-\psi)$ is locally nilpotent). Namely, a 
$(U_\hbar(\Ln_-),\psi)$-integrable $U_\hbar$-module $L$ goes to the 
$Z(U_\hbar)$-module $Hom_{U_\hbar}(\CW_\hbar^-,L)$. Conversely, 
a $Z(U_\hbar)$-module $R$ goes to $\CW_\hbar^-\otimes_{Z(U_\hbar)}R$.
In particular, $\CW_\hbar^-$ goes to the free module $Z(U_\hbar)$.

For an $\hbar$-Harish-Chandra bimodule $M$ we will construct a canonical
isomorphism $\kappa_\hbar(M)\simeq 
Hom_{U_\hbar}(\CW_\hbar^-,M\otimes_{U_\hbar}\CW_\hbar^-)$. In effect, 
$L\mapsto M\otimes_{U_\hbar}L$ is a right-exact endofunctor of the category of
$(U_\hbar(\Ln_-),\psi)$-integrable $U_\hbar$-modules. Under Kostant's 
equivalence, this endofunctor goes to the convolution with the 
$Z(U_\hbar)$-bimodule $X$ which corresponds by Kostant to our endofunctor
applied to $\CW_\hbar^-$. In other words, 
$X=Hom_{U_\hbar}(\CW_\hbar^-,M\otimes_{U_\hbar}\CW_\hbar^-)$.
We have a tautological isomorphism $X\otimes_{Z(U_\hbar)}\CW_\hbar^-\iso
M\otimes_{U_\hbar}\CW_\hbar^-$. This yields the desired isomorphism
$X\iso\kappa_\hbar(M)$. In particular, for a free $\hbar$-Harish-Chandra
bimodule $M=Fr(V)$, we obtain $\phi(V)\otimes_{Z(U_\hbar)}\CW_\hbar^-\iso
V\otimes_{\kk}\CW_\hbar^-$.

Now let us compute 
$(\CM_\hbar(-\rho)\otimes V\otimes \CW_\hbar^-)\otimes_{U_\hbar}
\kk[\hbar]=
(\CM_\hbar(-\rho)\otimes(V\otimes \CW_\hbar^-))\otimes_{U_\hbar}\kk[\hbar]\iso
(\CM_\hbar(-\rho)\otimes(\CW_\hbar^-\otimes_{Z(U_\hbar)}\phi(V)))
\otimes_{U_\hbar}\kk[\hbar]\iso(\pi,\Id,\Id)^*\phi(V)$. The last isomorphism
arises from $(\CM_\hbar(-\rho)
\otimes\CW_\hbar^-)\otimes_{U_\hbar}\kk[\hbar]\iso
U_\hbar(\Lt)=\O(\Lt^*)$, since 
$\CM_\hbar(-\rho)=U_\hbar\otimes_{U_\hbar(\Lb)}
\kk[\hbar][\Lt](-\rho)$, and $\CW_\hbar^-=U_\hbar\otimes_{U_\hbar(\Ln_-)}\psi$.

On the other hand, 
$(\CM_\hbar(-\rho)\otimes V\otimes \CW_\hbar^-)\otimes_{U_\hbar}
\kk[\hbar]=
((\CM_\hbar(-\rho)\otimes V)\otimes \CW_\hbar^-)
\otimes_{U_\hbar}\kk[\hbar]\iso
(\CM_\hbar(-\rho)\otimes V)\overset{L}{\otimes}_{U_\hbar(\Ln_-)_2}\psi\iso
\kk[\hbar][\Lt](-\rho)\overset{L}{\otimes}_{U_\hbar(\Lb)_1}(U_\hbar\otimes V)
\overset{L}{\otimes}_{U_\hbar(\Ln_-)_2}\psi=\varphi(V)$.

This completes the proof of the lemma.
\hfill $\Box$

\subsection{Canonical filtration on $\varphi(V)$}
\label{canon}
For $V\in Rep(\LG)$ we have $\varphi(V)=(\CM_\hbar(-\rho)\otimes V)
\otimes_{U_\hbar(\Ln_-)}\psi$. Note that $\CM_\hbar(-\rho)\otimes V$ has a
canonical filtration with associated graded 
$\bigoplus_\lambda\CM_\hbar(\lambda-\rho)\otimes\ _\lambda V$, where
$\lambda$ is a weight of $\Lt$, and $_\lambda V$ is the corresponding
weight space of $V$; furthermore,
$\CM_\hbar(\lambda-\rho)=
U_\hbar\otimes_{U_\hbar(\Lb)}\kk[\hbar][\Lt](\lambda-\rho)$,
and $\kk[\hbar][\Lt](\lambda-\rho)$ is a $U_\hbar(\Lb)$-module which factors
through the $U_\hbar(\Lt)=\kk[\hbar][\Lt]$-module where $t\in\Lt$ acts by
multiplication by $t+\hbar\lambda(t)-\hbar\rho(t)$.

It follows that $\varphi(V)$ has a canonical filtration with associated
graded $\bigoplus_\lambda
(\CM_\hbar(\lambda-\rho)\otimes_{U_\hbar(\Ln_-)}\psi)\otimes\ _\lambda V$.
Note that $\CM_\hbar(\lambda-\rho)\otimes_{U_\hbar(\Ln_-)}\psi$ is a
$\O(\Lt^*\times(\Lt^*/W)\times\Aone)$-module since
$\CM_\hbar(\lambda-\rho)$ is a $U_\hbar(\Lt)-Z(U_\hbar)$-bimodule.
To describe $\CM_\hbar(\lambda-\rho)\otimes_{U_\hbar(\Ln_-)}\psi$ as a coherent
sheaf on $\Lt^*\times(\Lt^*/W)\times\Aone$, we denote by 
$(\Id,\pi,\Id)$ the projection from $\Lt^*\times\Lt^*\times\Aone$ to
$\Lt^*\times(\Lt^*/W)\times\Aone$, and we denote by 
$\Gamma_\lambda\subset
\Lt^*\times\Lt^*\times\Aone$ the subscheme defined by the equations
$\Gamma_\lambda=\{(t_1,t_2,a):\ t_2=t_1+a\lambda\}$. Then 
$\CM_\hbar(\lambda-\rho)\otimes_{U_\hbar(\Ln_-)}\psi=
(\Id,\pi,\Id)_*\O(\Gamma_\lambda)$.

We have proved the following 

\begin{lem}
\label{Mit}
For $V\in Rep(\LG)$, the $\O(\Lt^*\times(\Lt^*/W)\times\Aone)$-module
$\varphi(V)$ has a canonical filtration with associated graded
$\bigoplus_\lambda
(\Id,\pi,\Id)_*\O(\Gamma_\lambda)\otimes\ _\lambda V$.
In particular, $\varphi(V)$ is flat as an $\O(\Lt^*\times\Aone)$-module.
\end{lem}
\hfill $\Box$

\subsection{Whittaker modules for Levi subalgebras}
Let $\LT\subset\LL\subset\LG$ be a Levi subgroup with the Lie algebra
$\Lt\subset\Ll\subset\Lg$. We denote by $\Lp_L$ (resp. $\Lp_L^-$) the 
parabolic subalgebra generated by $\Ll$ and the positive (resp. negative)
Borel subalgebra $\Lb$ (resp. $\Lb_-$). We denote by $\pi_L$ the 
projection from $\Lt^*/W_L$ to $\Lt^*/W$.

\begin{lem}
\label{Eugene}
For $V\in Rep(\LG)$, the $\O(\Lt^*/W_L\times\Lt^*/W\times\Aone)$-module
$(\pi_L,\Id,\Id)^*\phi(V)$ 
carries a canonical filtration $F^\bullet_L$ such that the associated graded is
equipped with a canonical isomorphism
$\ ^{\operatorname{alg}}\Xi_L:\ gr(\pi_L,\Id,\Id)^*\phi(V)\iso(\Id,\pi_L,\Id)_*\phi_L(V|_{\LL})$.
\end{lem}

\proof:
We have $\Ll=[\Ll,\Ll]\oplus\z_{\Ll}$ where $\z_{\Ll}$ stands for the
center of $\Ll$. We consider the nilpotent subalgebra $\Ln_-^L=\Ll\cap\Ln_-$,
and a nondegenerate homomorphism $\psi_L:\ U_\hbar(\Ln_-^L)\to\kk[\hbar]$
such that $\psi(f_\alpha)=1$ for any simple root $\alpha$ of $\Ll$.
We define the Whittaker $U_\hbar([\Ll,\Ll])$-module $\CW^-_L$ as
$U_\hbar([\Ll,\Ll])\otimes_{U_\hbar(\Ln_-^L)}\psi_L$. 
We define a free $U_\hbar(\z_{\Ll})=\kk[\hbar][\z_{\Ll}]$-module
$\z_{\Ll}(-\rho+\rho_L)$ as $\kk[\hbar][\z_{\Ll}]$ where $t\in\z_{\Ll}$
acts by multiplication by $t-\hbar(\rho-\rho_L)(t)$ (here $\rho_L$ is
the halfsum of positive roots of $\Ll$). 
We define a $U_\hbar(\Ll)$-module $\CW^-_L(-\rho+\rho_L)$ as
$\CW^-_L\otimes_{\kk[\hbar]}\z_{\Ll}(-\rho+\rho_L)$.
The projection $\Lp_L\to\Ll$ gives rise to the homomorphism
$U_\hbar(\Lp_L)\to U_\hbar(\Ll)$, and thus we can consider
$\CW^-_L(-\rho+\rho_L)$ as a $U_\hbar(\Lp_L)$-module.
Finally, we define the {\em Verma-Whittaker} $U_\hbar$-module
$\CM\CW_L(-\rho+\rho_L)$ as 
$U_\hbar\otimes_{U_\hbar(\Lp_L)}\CW^-_L(-\rho+\rho_L)$. Note that the
center $Z(U_\hbar(\Ll))=\O(\Lt^*/W_L\times\Aone)$ acts by endomorphisms
of $\CW^-_L(-\rho+\rho_L)$, and hence of $\CM\CW_L(-\rho+\rho_L)$.

We claim that for $V\in Rep(\LG)$, we have a canonical isomorphism
$$(\pi_L,\Id,\Id)^*\phi(V)[n_L]\cong
(\CM\CW_L(-\rho+\rho_L)\otimes 
V\otimes\CW^-_\hbar)\overset{L}{\otimes}_{U_\hbar}\kk[\hbar]$$ (the LHS
is homologically shifted to the degree $-n_L$, that is negative dimension
of $\Ln_-^L$).
In effect, arguing like in the proof of Lemma~\ref{mit}, we only have
to check that $(\CM\CW_L(-\rho+\rho_L)\otimes\CW^-_\hbar)
\overset{L}{\otimes}_{U_\hbar}
\kk[\hbar]\iso Z(U_\hbar(\Ll))[n_L]=\O(\Lt^*/W_L\times\Aone)[n_L]$.
To this end we note that $(\CM\CW_L(-\rho+\rho_L)\otimes\CW^-_\hbar)
\overset{L}{\otimes}_{U_\hbar}\kk[\hbar]\iso
\CW^-_L(-\rho+\rho_L)\overset{L}{\otimes}_{U_\hbar(\Lp_L)}
(\CW^-_\hbar|_{U_\hbar(\Lp_L)})=
\CW^-_L(-\rho+\rho_L)\overset{L}{\otimes}_{U_\hbar(\Lp_L)}
(U_\hbar(\Lp_L)\otimes_{U_\hbar(\Ln_-^L)}\psi_L)\iso
\CW^-_L(-\rho+\rho_L)\overset{L}{\otimes}_{U_\hbar(\Ln_-^L)}\psi_L\iso
Z(U_\hbar(\Ll))[n_L]$.

Moreover, it follows that for an $\LL$-module $W$ we have a canonical
isomorphism of $\O(\t/W_L\times\t/W\times\Aone)$-modules
$\left([U_\hbar\otimes_{U_\hbar(\Lp_L)}(\CW^-_L(-\rho+\rho_L)\otimes W)]
\otimes\CW^-_\hbar\right)\overset{L}{\otimes}_{U_\hbar}\kk[\hbar]\iso
(\Id,\pi_L,\Id)_*\phi_L(W)$. Now it remains to notice that for a $\LG$-module
$V$ the $U_\hbar$-module $\CM\CW_L(-\rho+\rho_L)\otimes V$ has a canonical
filtration with associated graded 
$U_\hbar\otimes_{U_\hbar(\Lp_L)}(\CW^-_L(-\rho+\rho_L)\otimes V|_{\LL})$.
This completes the proof of the lemma.
\hfill $\Box$

\subsection{Transitivity for a pair of Levi subgroups}
\label{algpair}
We have a canonical isomorphism
$\ ^{\operatorname{alg}}\Xi_L:\ gr(\pi_L,\Id,\Id)^*\phi(V)
\iso(\Id,\pi_L,\Id)_*\phi_L(V|_{\LL})$. 
In the RHS we have the restriction of 
$\O(\Lt^*/W_L\times\Lt^*/W_L\times\Aone)$-module $\phi_L(V|_{\LL})$ to
$\O(\Lt^*/W_L\times\Lt^*/W\times\Aone)$. To save a bit of notation in what 
follows we will write simply
$\ ^{\operatorname{alg}}\Xi_L:\ gr(\pi_L,\Id,\Id)^*\phi(V)\iso\phi_L(V|_{\LL})$.

If $\LT\subset\LL'\subset\LL$
is another Levi subgroup, then we denote by $\pi_{L'}^L$ the projection
from $\Lt^*/W_{L'}$ to $\Lt^*/W_L$. Note that the filtration $F_{L'}^\bullet$
on $(\pi_{L'},\Id,\Id)^*\phi(V)=
(\pi_{L'}^L,\Id,\Id)^*(\pi_L,\Id,\Id)^*\phi(V)$
is a refinement of the filtration $(\pi_{L'}^L,\Id,\Id)^*F_L^\bullet$,
and hence induces a canonical filtration $F_{L'}^{L\bullet}$ on
$(\pi_{L'}^L,\Id,\Id)^*gr_{F_L^\bullet}(\pi_L,\Id,\Id)^*
\phi(V)$. The isomorphism
$(\pi_{L'}^L,\Id,\Id)^*\ ^{\operatorname{alg}}\Xi_L$ carries the filtration $F_{L'}^{L\bullet}$
to the filtration $F_{L'}^\bullet$ on
$(\pi_{L'}^L,\Id,\Id)^*\phi_L(V|_{\LL})$.
We have a canonical isomorphism
$\ ^{\operatorname{alg}}\Xi_{L'}^L:\ gr_{F_{L'}^\bullet}(\pi_{L'}^L,\Id,\Id)^*
\phi_L(V|_{\LL})\iso\phi_{L'}(V|_{\LL'})$.
We consider the composition
\begin{equation}
\begin{CD} 
gr_{F_{L'}^\bullet}(\pi_{L'},\Id,\Id)^*\phi(V)
@>gr_{F_{L'}^{L\bullet}}(\pi_{L'}^L,\Id,\Id)^*\ ^{\operatorname{alg}}\Xi_L>>
%\end{CD}
%\end{equation}
%$$\begin{CD}
gr_{F_{L'}^\bullet}(\pi_{L'}^L,\Id,\Id)^*
\phi_L(V|_{\LL})
@>\ ^{\operatorname{alg}}\Xi_{L'}^L>>
\phi_{L'}(V|_{\LL'})
\end{CD}
\end{equation}
Then we have
\begin{lem}
\label{ap}
$^{\operatorname{alg}}\Xi_{L'}^L\circ gr_{F_{L'}^{L\bullet}}
(\pi_{L'}^L,\Id,\Id)^*\ ^{\operatorname{alg}}\Xi_L=\ 
^{\operatorname{alg}}\Xi_{L'}$.
\end{lem}
\hfill $\Box$

\subsection{Tensor structure on Kostant functor}
\label{tensphi}
Recall that according to the proof of Lemma~\ref{mit}, for $V\in Rep(\LG)$,
we have a canonical isomorphism $\phi(V)\otimes_{\O(\Lt^*/W\times\Aone)}
\CW^-_\hbar\iso V\otimes\CW^-_\hbar$. Thus, for $V_1,V_2\in Rep(\LG)$,
we have $\phi(V_1)\otimes_{\O(\Lt^*/W\times\Aone)}\phi(V_2)
\otimes_{\O(\Lt^*/W\times\Aone)}\CW^-_\hbar\iso
\phi(V_1)\otimes_{\O(\Lt^*/W\times\Aone)}(V_2\otimes\CW^-_\hbar)\iso
(V_1\otimes V_2)\otimes\CW^-_\hbar\stackrel{\sim}{\longleftarrow}
\phi(V_1\otimes V_2)\otimes_{\O(\Lt^*/W\times\Aone)}\CW^-_\hbar$.
Composing (and inverting) these isomorphisms we obtain
$\phi(V_1\otimes V_2)\otimes_{\O(\Lt^*/W\times\Aone)}\CW^-_\hbar\iso
\phi(V_1)\otimes_{\O(\Lt^*/W\times\Aone)}\phi(V_2)
\otimes_{\O(\Lt^*/W\times\Aone)}\CW^-_\hbar$, and thus
$$\ ^{\operatorname{alg}}\omega_{V_1,V_2}:\ \phi(V_1\otimes V_2)\iso
\phi(V_1)\otimes_{\O(\Lt^*/W\times\Aone)}\phi(V_2)=:\phi(V_1)\star\phi(V_2).$$
According to Lemma~\ref{Mit} (cf. also Lemma~\ref{Eugene}),
we have a canonical isomorphism 
$^{\operatorname{alg}}\Xi_V=\ 
^{\operatorname{alg}}\Xi_{T,V}:\ gr(\pi,\Id,\Id)^*\phi(V)\iso
(\Id,\pi,\Id)_*\phi_T(V|_{\LT})$. In the RHS we have the restriction of
$\O(\Lt^*\times\Lt^*\times\Aone)$-module $\phi_T(V|_{\LT})$ to
$\O(\Lt^*\times(\Lt^*/W)\times\Aone)$. To save a bit of notation in what
follows we will write simply 
$^{\operatorname{alg}}\Xi_V:\ gr(\pi,\Id,\Id)^*\phi(V)\iso
\phi_T(V|_{\LT})$. It follows that after tensoring with
$\kk(\Lt^*\times\Aone)$ (over the first and third factors in
$\O(\Lt^*/W\times\Lt^*/W\times\Aone)$) we have a canonical isomorphism
$$^{\operatorname{alg}}\Xi_V^\gen:\ 
\phi(V)\otimes_{\O(\Lt^*/W\times\Aone)}\kk(\Lt^*\times\Aone)=
gr(\pi,\Id,\Id)^*\phi(V)
\otimes_{\O(\Lt^*\times\Aone)}\kk(\Lt^*\times\Aone)\iso$$
$$\iso\phi_T(V|_{\LT})\otimes_{\O(\Lt^*\times\Aone)}\kk(\Lt^*\times\Aone)=
\bigoplus_\lambda\left(\O(\Gamma_\lambda)
\otimes_{\O(\Lt^*\times\Aone)}\kk(\Lt^*\times\Aone)\right)\otimes\ 
_\lambda V$$
Now we have a canonical isomorphism 
$\O(\Gamma_\mu)\star\O(\Gamma_\nu):=\O(\Gamma_\mu)
\otimes_{\O(\Lt^*\times\Aone)}\O(\Gamma_\nu)=\O(\Gamma_{\mu+\nu})$.
Hence we get a canonical isomorphism 
$^{\operatorname{alg}}\Xi_{V_1}^\gen\star\ 
^{\operatorname{alg}}\Xi_{V_2}^\gen:$ 
$$(\phi(V_1)\star\phi(V_2)) 
\otimes_{\O(\Lt^*/W\times\Aone)}\kk(\Lt^*\times\Aone)\iso
\bigoplus_{\mu+\nu=\lambda}\left(\O(\Gamma_\lambda)
\otimes_{\O(\Lt^*\times\Aone)}\kk(\Lt^*\times\Aone)\right)\otimes\
_\mu V_1\otimes\ _\nu V_2=$$
$$=\bigoplus_\lambda\left(\O(\Gamma_\lambda)
\otimes_{\O(\Lt^*\times\Aone)}\kk(\Lt^*\times\Aone)\right)\otimes\
_\lambda(V_1\otimes V_2)$$
We want to compare it with
$$\ ^{\operatorname{alg}}\Xi_{V_1\otimes V_2}^\gen:\
\phi(V_1\otimes V_2)\otimes_{\O(\Lt^*/W\times\Aone)}\kk(\Lt^*\times\Aone)
\iso\bigoplus_\lambda\left(\O(\Gamma_\lambda)
\otimes_{\O(\Lt^*\times\Aone)}\kk(\Lt^*\times\Aone)\right)\otimes\
_\lambda(V_1\otimes V_2)$$

\begin{prop}
\label{tensalg}
$\ ^{\operatorname{alg}}\Xi_{V_1\otimes V_2}^\gen=
(\ ^{\operatorname{alg}}\Xi_{V_1}^\gen\star\ 
^{\operatorname{alg}}\Xi_{V_2}^\gen)
\circ\ ^{\operatorname{alg}}\omega_{V_1,V_2}$.
\end{prop}

\proof: We consider the generic universal Verma module 
$\CMM=U_\hbar\otimes_{U_\hbar(\Lb)}\kk(\Lt^*\times\Aone)(-\rho)$, and
$\kk(\Lt^*\times\Aone)(-\rho)$ is a $U_\hbar(\Lb)$-module which factors
through the $U_\hbar(\Lt)=\O(\Lt^*\times\Aone)$-module where $t\in\Lt$ 
acts by multiplication by $t-\hbar\rho(t)$. It is well known that
$\End_{U_\hbar}(\CMM)=\kk(\Lt^*\times\Aone)$, and the category $\CB$ 
of $U_\hbar(\Ln)$-integrable $U_\hbar\otimes\kk(\Lt^*\times\Aone)$-modules is 
semisimple, and any simple object is isomorphic to $\CMM$. In particular,
$\CB$ is equivalent to the category of 
$\kk(\Lt^*\times\Aone)$-modules. 

For $V\in Rep(\LG)$ we put $\varphi^{\operatorname{gen}}(V):=
\varphi(V)\otimes_{\O(\Lt^*\times\Aone)}\kk(\Lt^*\times\Aone)=
gr\varphi(V)\otimes_{\O(\Lt^*\times\Aone)}\kk(\Lt^*\times\Aone)$.
This is the restriction of a $\kk(\Lt^*\times\Aone)\otimes\O(\Lt^*)$-module
to $\kk(\Lt^*\times\Aone)\otimes\O(\Lt^*/W)$, but we will view it as a
$\kk(\Lt^*\times\Aone)\otimes\O(\Lt^*)$-module.

Arguing like in the proof of Lemma~\ref{mit}, we
obtain a canonical isomorphism 
$\varphi^{\operatorname{gen}}(V)\otimes_{\kk(\Lt^*\times\Aone)}\CMM\iso
V\otimes\CMM$. This gives rise to the tensor structure on the functor
$\varphi^{\operatorname{gen}}$:
$$\varphi^{\operatorname{gen}}(V_1\otimes V_2)\iso
\varphi^{\operatorname{gen}}(V_1)\star\varphi^{\operatorname{gen}}(V_2).$$
Clearly, the identification 
$^{\operatorname{alg}}\Xi_V^\gen:\ \varphi^\gen(V)\iso
\bigoplus_\lambda\left(\O(\Gamma_\lambda)
\otimes_{\O(\Lt^*\times\Aone)}\kk(\Lt^*\times\Aone)\right)\otimes\ 
_\lambda V$ commutes with the obvious tensor structure in the RHS.

On the other hand, arguing like in the proof of Lemma~\ref{mit}, we
obtain a canonical isomorphism $\varphi^\gen(V)\iso
(\CMM\otimes V\otimes \CW^-_\hbar)\otimes_{U_\hbar}\kk[\hbar]$ which 
implies that the tensor structures on $\phi$ and $\varphi^\gen$ are
compatible as well. This completes the proof of the proposition. \hfill $\Box$

\subsection{Quasiclassical limit of $\phi(V)$}
\label{secr}
For $V\in Rep(\LG)$, Lemma~\ref{Mit} implies that the $\O(\Lt^*/W)$-bimodule 
$\phi(V)$ is supported at the diagonal $\Delta\subset\Lt^*/W\times\Lt^*/W$.
It follows that the action of $\O(\Lt^*/W\times\Lt^*/W\times\Aone)$ on 
$\phi(V)$ actually extends to the action of
$\O(N_{\Lt^*/W\times\Lt^*/W}\Delta)$. As we know from~\ref{quasic},
$\O(N_{\Lt^*/W\times\Lt^*/W}\Delta)/\hbar\simeq\O(\bT(\Lt^*/W))$ is the 
universal centralizer.

\begin{lem}
\label{MIt}
For $V\in Rep(\LG)$, the $\O(\bT(\Lt^*/W))$-module $\phi(V)|_{\hbar=0}$
is canonically isomorphic to the module $\O(\Sigma)\otimes V$ over the
universal centralizer.
\end{lem}

\proof: Consider the $\O(\Lg^*)^{\LG}$-module $Pol(\Lg^*,\Lg)^{\LG}$ of
$\LG$-invariant polynomial maps from $\Lg^*$ to $\Lg$; this is a vector
bundle over $Spec\O(\Lg^*)^{\LG}=\Sigma$. Given a polynomial 
$P\in\O(\Lg^*)^{\LG}$, its differential $dP$ defines a section of this
vector bundle. 
For a central element $z\in Z(U(\Lg))=\O(\Lg^*)^{\LG}$ we denote the
corresponding section by $\sigma_z$.
If $z$ runs through a set of generators of $Z(U(\Lg))$, the corresponding
sections $\sigma_z$ form a basis of the universal centralizer bundle,
and identify it with the cotangent bundle $\bT^*(\Sigma)$.

Thus it suffices to check the following statement about the free
$\hbar$-Harish-Chandra bimodule $U_\hbar\otimes V$. Let $z^{(1)}$
(resp. $z^{(2)}$) stand for the left (resp. right) action of $z$ in
$U_\hbar\otimes V$. Then the action of 
$\frac{z^{(1)}-z^{(2)}}{\hbar}|_{\hbar=0}$ on $(U_\hbar\otimes V)|_{\hbar=0}=
\O(\Lg^*)\otimes V$ coincides with the action of $\sigma_z\in
\O(\Lg^*)\otimes\Lg$.

In effect, if $v\in V$, and $z=\sum_{\underline{i}} 
a_{\underline{i}}x_{i_1}\ldots x_{i_k}$ where $x_{i_l}\in\Lg$, 
and $\tilde{y}\in U_\hbar$ is a lift of $y\in\O(\Lg^*)=U_\hbar|_{\hbar=0}$,
then $\frac{z(\tilde{y}\otimes v)-(\tilde{y}\otimes v)z}{\hbar}|_{\hbar=0}=
\sum_{i_l\in\underline{i}}
a_{\underline{i}}x_{i_1}\ldots\widehat{x_{i_l}}\ldots 
x_{i_k}y\otimes x_{i_l}(v)=\sigma_z(y\otimes v)$.

The lemma is proved. 
\hfill $\Box$

\section{Rank 1}

\subsection{Equivariant cohomology for $G=PGL(2)$}
\label{before}
Let us describe 
$H^\bullet_{\GO\rtimes\Gm}(\Gr_G,S(V))$ in the case $G=PGL(2)$. 
Then $\LG$ is isomorphic to $SL(2)$.
For a positive integer $n$ let $V_n$ be the $n+1$-dimensional irreducible
$SL(2)$-module. Let $\Gr_n$ be the closure of the $n$-dimensional $\GO$-orbit
in $\Gr_G$. It is known that $\Gr_n$ is rationally smooth, so we are 
interested in $H^\bullet_{\GO\rtimes\Gm}(\Gr_n)$ as a module over
$H^\bullet_{\GO\rtimes\Gm}(\Gr)$. The non-equivariant cohomology
$H^\bullet(\Gr_n)$ is a cyclic $H^\bullet(\Gr)$-module (see e.g.~\cite{G1}).
Hence, by graded Nakayama lemma, $H^\bullet_{\GO\rtimes\Gm}(\Gr_n)$ is a cyclic
module over $H^\bullet_{\GO\rtimes\Gm}(\Gr)$. Recall that
$H^\bullet_{\GO\rtimes\Gm}(\Gr)\cong\O(N_{\t/W\times\t/W}\Delta)\oplus
\O(N_{\t/W\times\t/W}\Delta)$, and so $H^\bullet_{\GO\rtimes\Gm}(\Gr_n)$
is the structure sheaf of a subscheme $A_n\subset N_{\t/W\times\t/W}\Delta$
of a copy of $N_{\t/W\times\t/W}\Delta$ specified by the parity of $n$.

Now we describe the subscheme $A_n$. Let $P_n=\{n\omega,(n-2)\omega,\ldots,
(2-n)\omega,-n\omega\}$ be the set of weights of $\LG$-module $V_n$.
We have $P_n\subset\t=\Lt^*$. For $i=-n,-n+2,\ldots,n-2,n$,
let $\Gamma_i\subset\t\times\t\times\Aone$
be a subscheme defined by the equations $\Gamma_i=\{(x_1,x_2,a):\
x_2=x_1+ia\omega\}$. Let $\Gamma(n)$ stand for the subscheme defined by
the product of the above equations (over $i=-n,-n+2,\ldots,n-2,n$).
Recall that $\pi$ stands for the projection $\t\to\t/W$,
and consider the subscheme 
$(\pi,\pi,\Id)(\Gamma(n))\subset\t/W\times\t/W\times\Aone$. 
Finally, we can formulate

\begin{lem}
\label{top}
$A_n$ is the strict transform of $(\pi,\pi,\Id)(\Gamma(n))$ in
$N_{\t/W\times\t/W}\Delta$.
\end{lem}

\proof: 
Since $H^\bullet_{\GO\rtimes\Gm}$ is a flat $\O(\t/W\times\Aone)$-module,
it suffices to identify $A_n$ with $(\pi,\pi,\Id)(\Gamma(n))$ generically
over $\t/W\times\Aone$, or else to identify $(\pi,\Id,\Id)^{-1}(A_n)$
with $(\Id,\pi,\Id)(\Gamma(n))$ generically over $\t\times\Aone$.
This was done in Lemma~\ref{Harvard}.
\hfill $\Box$

\subsection{Generic splitting of the canonical filtration on equivariant
cohomology for $G=PGL(2)$}
\label{?}
Recall the canonical filtration on
$H^\bullet_{T\times\Gm}(\Gr_n)$ (see~\ref{can}). 
We will compare the identification
(see Lemma~\ref{Harvard}) of the associated graded with 
$\bigoplus_{i=-n}^n(\Id,\pi,\Id)_*\O(\Gamma_i)\otimes\ _iV_n$ with the
identification (see Lemma~\ref{top}) 
$H^\bullet_{T\times\Gm}(\Gr_n)\cong\O\left((\pi,\Id,\Id)^{-1}(A_n)\right)$.
To this end we recall some basic facts about the cohomology of $\Gr_n$.
For $i=-n,-n+2,\ldots,n-2,n$, let $v_i\in H^{i+n}(\Gr_n)$ stand for 
the (Poincar{\'e} dual of the) 
fundamental class of $\Gr_n\cap\ofT_i$ (an irreducible subvariety of
$\Gr_n$ of dimension $\frac{n-i}{2}$). 
The action of $e,h,f\in\mathfrak{sl}_2$ on
$H^\bullet(\Gr_n)$ in this basis is given by
$$hv_i=iv_i,\ ev_{i-2}=\frac{n+i}{2}v_i,\ fv_{i+2}=\frac{n-i}{2}v_i$$
(recall that $e$ is defined as the multiplication by the first Chern class
of the determinant line bundle). 

The canonical filtration $0=F^{n+2}\subset F^n\subset F^{n-2}\subset
\ldots\subset F^{2-n}\subset F^{-n}=H^\bullet_{T\times\Gm}(\Gr_n)$ is given
by $F^i=\operatorname{Im}\left(r_i:\ H^\bullet_{\ofT_i,T\times\Gm}(\Gr_n)\to
H^\bullet_{T\times\Gm}(\Gr_n)\right)$.
On the other hand, the proof of Lemma~\ref{top} shows that 
$H^\bullet_{T\times\Gm}(\Gr_n)$ is generated by the (Poincar{\'e} dual of the)
fundamental class $\tilde{v}_{-n}$ of $\Gr_n=\Gr_n\cap\ofT_{-n}$, i.e.\ by the unit cohomology class~1.
Recall that $\fj_i$ stands for the embedding of the $T$-fixed point
$i$ into $\fT_i$, while $\imath_i$ stands for the locally closed 
embedding of $\fT_i$ into $\Gr$. Let us denote by $\bar{\imath}_i$
the closed embedding of $\ofT_i$ into $\Gr$.
Let $\bar s_i$ stand for the unit element in the equivariant hyperbolic stalk $\fj_i^*\imath_i^!\underline\BC{}_{\Gr_n}$
(it lives in cohomological degree $i+n$). Then under the isomorphism of the Localization theorem for $T\times{\mathbb G}_m$-equivariant
cohomology, we have 
\begin{equation}
\label{Idiot}
1=\tilde v_{-n}=\sum_{i=-n}^n(x+(i-1)\hbar)^{-1}(x+(i-2)\hbar)^{-1}\ldots(x+\frac{i-n}{2}\hbar)^{-1}\bar s_i.
\end{equation}
Indeed,
%The image of $\bar{\imath}{}^*_i\tilde{v}_{-n}$ in the nonequivariant
%cohomology $H^\bullet(\Gr_n\cap\ofT_i)$ is the fundamental class $v_i$ of
%$\Gr_n\cap\ofT_i$. We can further restrict it to $\Gr_n\cap\fT_i$,
%and then to $i$ to get $\fj_i^*\imath_i^*\tilde{v}_{-n}$.
%To compare it with $\fj_i^*\imath_i^!\tilde{v}_{-n}$
we consider a transversal slice $\Gr_n\cap\fS_i$ to $\Gr_n\cap\ofT_i$ in $\Gr_n$
where $\fS_i$ is the $N({\mathbf F})$-orbit through the point $i$.
It is known that $\Gr_n\cap\fS_i$ is isomorphic to a vector space 
$\BA^{\frac{n+i}{2}}$ with the origin at $i$, and the action of $T\times\Gm$ is
linear with weights $x+(i-1)\hbar,x+(i-2)\hbar,\ldots,x+\frac{i-n}{2}\hbar$.
%It follows that $\imath_i^!\tilde{v}_{-n}=(x+(i-1)\hbar)(x+(i-2)\hbar)\ldots
%(x+\frac{i-n}{2}\hbar)\imath_i^*\tilde{v}_{-n}$. We conclude that the generator $\tilde{v}_i$ 
%of $F^i$ whose class in the nonequivariant cohomology $H^\bullet(\Gr_n)$ is equal to $v_i$ is given by

\subsection{Kostant functor for $\LG=SL(2)$}
Now we consider the group $\LG=SL(2)$ with the Lie algebra 
$\Lg=\mathfrak{sl}_2$ and Cartan subalgebra $\Lt\subset\Lg$.

\begin{lem}
\label{alg}
The $\O(\Lt^*/W\times\Lt^*/W\times\Aone)$-module
$\phi(V_n)=\kappah(U_\hbar\otimes V_n)$ is isomorphic to $\O(A_n)$.
\end{lem}

\proof: According to Lemma~\ref{MIt}, 
the restriction of $\phi(V_n)$ to $\hbar=0$ is isomorphic to the
$\O(\bT(\Lt^*/W))$-module $V_n\otimes\O(\Lt^*/W)$, that is $V_n$ viewed
as a module over the universal centralizer. Further restricting it to 
$0\in\Lt^*/W$ we obtain $V_n$ viewed as a module over the centralizer of the
regular nilpotent $e\in\mathfrak{sl}_2$. Clearly, $V_n$ is a cyclic
$\kk[e]$-module. By the graded Nakayama Lemma, $\phi(V_n)$ is a cyclic
$\O(N_{\Lt^*/W\times\Lt^*/W}\Delta)$-module as well, hence $\phi(V_n)$ is
isomorphic to the structure sheaf of a subscheme 
$B_n\subset N_{\Lt^*/W\times\Lt^*/W}\Delta$. 
We have to check $B_n=A_n$. 

According to Lemmas~\ref{mit},~\ref{Mit}, $\phi(V_n)$ is a flat
$\O(\Lt^*/W\times\Aone)$-module, so it suffices to identify $B_n$ and
$A_n$ generically over $\Lt^*/W\times\Aone$. Moreover, it suffices to 
identify $(\pi,\pi,\Id)^{-1}(B_n)$ with 
$\Gamma(n)\subset\Lt^*\times\Lt^*\times\Aone$.
This was done in Lemma~\ref{Mit}.
This completes the proof of the lemma.
\hfill $\Box$

\subsection{Generic splitting of the canonical filtration on Kostant functor
for $\LG=SL(2)$} Recall the canonical filtration on
$\varphi(V_n)$ (see~\ref{canon}). We will compare the identification
(see Lemma~\ref{Mit}) of the associated graded with 
$\bigoplus_{i=-n}^n(Id,\pi,\Id)_*\O(\Gamma_i)\otimes\ _iV_n$ with the
identification (see Lemma~\ref{alg}) 
$\varphi(V_n)\cong\O\left((\pi,\Id,\Id)^{-1}(A_n)\right)$.
To this end we recall some basic facts about 
$U_\hbar(\mathfrak{sl}_2)$-modules. First, $V_n$ is a free $\kk[\hbar]$-module
with a basis $\{v_n,v_{n-2},\ldots,v_{2-n},v_{-n}\}$. The action of
$e,h,f\in\mathfrak{sl}_2$ is given by
$$hv_i=i\hbar v_i,\ ev_{i-2}=\frac{n+i}{2}\hbar v_i,\
fv_{i+2}=\frac{n-i}{2}\hbar v_i$$
Second, $\CM_\hbar(-1)$ is a free $\kk[\hbar,x]$-module with a basis
$\{m_{-1},m_{-3},m_{-5},\ldots\}$. The action of $e,h,f$ is given by
$$hm_i=(x+i\hbar)m_i,\ 
em_i=\frac{-i-1}{2}\hbar(x-\frac{i+1}{2}\hbar)m_{i+2},\ 
fm_i=m_{i-2}$$

We have a canonical filtration 
$0=F^{n+2}\subset F^n\subset F^{n-2}\subset\ldots\subset
F^{2-n}\subset F^{-n}=
\CM_\hbar(-1)\otimes_{\kk[\hbar]}V_n$ by $U_\hbar$-submodules
such that $F^i/F^{i+2}=\CM_\hbar(i-1)\otimes\ _iV_n$ 
(notations of~\ref{canon}).
Recall that $_iV_n$ is spanned by $v_i$. There is a unique vector
$s_i\in(\CM_\hbar(-1)\otimes V_n)\otimes_{\kk[\hbar,x]}\kk(\hbar,x)$
such that $es_i=0$, and $s_i\equiv m_{-1}\otimes v_i$ modulo
$U_\hbar\langle s_{i+2},s_{i+4},\ldots,s_n\rangle$. 
Then $F^i=U_\hbar\langle s_i,s_{i+2},\ldots,s_n\rangle\cap
(\CM_\hbar(-1)\otimes V_n)$. The image $\bar{s}_i$ of this vector in
the $\psi$-coinvariants 
$U_\hbar\langle s_i\rangle\otimes_{U_\hbar(\Ln_-)}\psi=
\CM_\hbar(i-1)\otimes_{U_\hbar(\Ln_-)}\psi$ is the generator of 
$(\Id,\pi,\Id)_*\O(\Gamma_i)\otimes\ _iV_n$. Note that the space
of $\psi$-coinvariants is just the quotient modulo the image of $f-1$.

On the other hand, the proof of Lemma~\ref{alg} shows that
$\varphi(V_n)=
(\CM_\hbar(-1)\otimes V_n)\otimes_{U_\hbar(\Ln_-)}\psi$ is generated by
the image $\overline{m_{-1}\otimes v_{-n}}$ of $m_{-1}\otimes v_{-n}$ in the space
of $\psi$-coinvariants. Thus we have to express 
$\overline{m_{-1}\otimes v_{-n}}$
in terms of $\bar{s}_{-n},\ldots,\bar{s}_n$.

\begin{lem}
\label{idiot}
$$\overline{m_{-1}\otimes v_{-n}}=\sum_{i=-n}^n
(x+(i-1)\hbar)^{-1}(x+(i-2)\hbar)^{-1}
\ldots\left(x+\frac{i-n}{2}\hbar\right)^{-1}
\bar{s}_i$$
\end{lem}

\proof: Recall that the space
of $\psi$-coinvariants is just the quotient modulo the image of $f-1$.
This means we have to prove the following equality in
$(\CM_\hbar(-1)\otimes V_n)\otimes_{\kk[\hbar,x]}\kk(\hbar,x)$:
$$m_{-1}\otimes v_{-n}=\sum_{i=-n}^n%\hbar^{\frac{-n-i}{2}}
(x+(i-1)\hbar)^{-1}(x+(i-2)\hbar)^{-1}
\ldots\left(x+\frac{i-n}{2}\hbar\right)^{-1}
f^{\frac{n+i}{2}}s_i$$
For $l=0,1,\ldots,n$ we introduce a new vector $s_i^l$ such that
$es_i^l=0$, and $s_i^l\equiv e^l(m_{-1}\otimes v_{i-2l})$ modulo
$U_\hbar\langle s_{i+2},s_{i+4},\ldots,s_{n}\rangle$ (evidently,
$s_i^l$ is proportional to $s_i=s_i^0$ 
when $i-2l\geq-n$; otherwise, $s_i^l$ is 
not defined). Consider the following collection of equalities:
$$e^l(m_{-1}\otimes v_{-n})=\sum_{i=-n+2l}^n\hbar^l%\hbar^{\frac{2l-n-i}{2}}
(x+(i-1)\hbar)^{-1}(x+(i-2)\hbar)^{-1}
\ldots\left(x+\frac{i-n+2l}{2}\hbar\right)^{-1}
f^{\frac{n-2l+i}{2}}s_i^l$$
Then the $n$-th equality is obvious, while the $l+1$-st equality is equivalent
to the $l$-th one by applying $e$ to both sides. Thus the desired
equality (for $l=0$) follows by descending induction in $l$.
\hfill $\Box$

\section{Topology vs Algebra}

\subsection{Comparison of Kostant functor with equivariant cohomology}

Recall that $\t$ is identified with $\Lt^*$.

\begin{thm}
\label{compaThm}
a) For $V\in Rep(\LG)$ there is a unique isomorphism of
$\O(\t/W\times\t/W\times\Aone)$-modules
$\eta_V:\ H^\bullet_{\GO\rtimes\Gm}(\Gr_G,S(V))\simeq\phi(V)$ 
such that $(\pi,\Id,\Id)^*\eta_V$ preserves the
canonical filtrations and induces the identity isomorphism of the
associated graded 
$gr(\pi,\Id,\Id)^*H^\bullet_{\GO\rtimes\Gm}(\Gr_G,S(V))=
(\Id,\pi,\Id)_*\bigoplus_\lambda\O(\Gamma_\lambda)\otimes\ _\lambda V=
gr(\pi,\Id,\Id)^*\phi(V)$. 

b) For $V_1,V_2\in Rep(\LG)$ the composition
$$\begin{CD} 
H^\bullet_{\GO\rtimes\Gm}(\Gr_G,S(V_1))
\star H^\bullet_{\GO\rtimes\Gm}(\Gr_G,S(V_2))
@>\ ^{\operatorname{top}}\omega_{V_1,V_2}^{-1}>>
H^\bullet_{\GO\rtimes\Gm}(\Gr_G,S(V_1\otimes V_2))
\end{CD}$$
$$\begin{CD}
@>\eta_{V_1\otimes V_2}>>
\phi(V_1\otimes V_2)
@>\ ^{\operatorname{alg}}\omega_{V_1,V_2}>>
\phi(V_1)\star\phi(V_2)
\end{CD}$$ equals $\eta_{V_1}*\eta_{V_2}$ (notations  
of~\ref{tenscoh} and~\ref{tensphi}).
\end{thm}

\proof: a) We have $$(\pi,\Id,\Id)^*H^\bullet_{\GO\rtimes\Gm}(\Gr_G,S(V))
\otimes_{\O(\t\times\Aone)}\kk(\t\times\Aone)=$$
$$gr(\pi,\Id,\Id)^*H^\bullet_{\GO\rtimes\Gm}(\Gr_G,S(V))
\otimes_{\O(\t\times\Aone)}\kk(\t\times\Aone)=$$
$$\left(\bigoplus_\lambda(\Id,\pi,\Id)_*\O(\Gamma_\lambda)\otimes\ _\lambda V
\right)\otimes_{\O(\t\times\Aone)}\kk(\t\times\Aone)=$$
$$gr(\pi,\Id,\Id)^*\phi(V)\otimes_{\O(\t\times\Aone)}\kk(\t\times\Aone)=
(\pi,\Id,\Id)^*\phi(V)\otimes_{\O(\t\times\Aone)}\kk(\t\times\Aone).$$
Thus we have to identify $(\pi,\Id,\Id)^*H^\bullet_{\GO\rtimes\Gm}(\Gr_G,S(V))$
with the natural $W$-action on it, and $(\pi,\Id,\Id)^*\phi(V)$ with the
natural $W$-action on it, as two $\O(\t\times\t/W\times\Aone)$-submodules
of $\left(\bigoplus_\lambda(\Id,\pi,\Id)_*\O(\Gamma_\lambda)\otimes\ _\lambda V
\right)\otimes_{\O(\t\times\Aone)}\kk(\t\times\Aone)$. First we will
show that the $W$-action on
$\left(\bigoplus_\lambda(\Id,\pi,\Id)_*\O(\Gamma_\lambda)\otimes\ _\lambda V
\right)\otimes_{\O(\t\times\Aone)}\kk(\t\times\Aone)$ 
arising from its identification with
$(\pi,\Id,\Id)^*H^\bullet_{\GO\rtimes\Gm}(\Gr_G,S(V))
\otimes_{\O(\t\times\Aone)}\kk(\t\times\Aone)$ coincides with the $W$-action
arising from the identification with 
$(\pi,\Id,\Id)^*\phi(V)\otimes_{\O(\t\times\Aone)}\kk(\t\times\Aone)$.
 
Let $\alpha$ be a simple root, $s_\alpha\in W$ the corresponding
simple reflection, and $(e_\alpha,h_\alpha,f_\alpha)$ the corresponding
$\mathfrak{sl}_2$-triple in $\Lg$. Let $v\in\ _\lambda V$ be a vector such that
$f_\alpha v=0$. Then $h_\alpha v=\lambda(h_\alpha)v$, and $\lambda(h_\alpha)$
is a nonpositive integer. We consider the following vector in
$\left(\bigoplus_\lambda(\Id,\pi,\Id)_*\O(\Gamma_\lambda)\otimes\ _\lambda V
\right)\otimes_{\O(\t\times\Aone)}\kk(\t\times\Aone)\ni
p_\alpha(v):=$
$$\sum_{k=0}^{-\lambda(h_\alpha)}
(h_\alpha+(2k-1+\lambda(h_\alpha))\hbar)^{-1}
(h_\alpha+(2k-2+\lambda(h_\alpha))\hbar)^{-1}\ldots
(h_\alpha+(k+\lambda(h_\alpha))\hbar)^{-1}
1_{\lambda+k\alpha}\otimes e_\alpha^kv$$
where $1_\mu$ is the constant function 1 on $(\Id,\pi,\Id)\Gamma_\mu$.

It follows from the computation in rank 1 (Lemma~\ref{idiot}),
and the transitivity equation~(\ref{ap}) of subsection~\ref{algpair}
(applied to the case where $\LL$ is a subminimal Levi containing just one
positive root $\alpha$, and $\LL'=\LT$) 
that $p_\alpha(v)$ is $s_\alpha$-invariant
(with respect to the action arising from the identification with
$(\pi,\Id,\Id)^*\phi(V)\otimes_{\O(\t\times\Aone)}\kk(\t\times\Aone)$),
and $M_\alpha:=
\sum_\lambda^{v\in\ _\lambda V^{f_\alpha}}\O(\t\times\Aone)p_\alpha(v)$
contains $(\pi,\Id,\Id)^*\phi(V)$.

On the other hand, it follows from the computation in rank 1 
(equation~(\ref{Idiot}) of subsection~\ref{?}),
and the transitivity equation~(\ref{tp}) of subsection~\ref{toppair}
(applied to the case where $\LL$ is a subminimal Levi containing just one
positive root $\alpha$, and $\LL'=\LT$) 
that $p_\alpha(v)$ is $s_\alpha$-invariant
(with respect to the action arising from the identification with
$(\pi,\Id,\Id)^*H^\bullet_{\GO\rtimes\Gm}(\Gr_G,S(V))
\otimes_{\O(\t\times\Aone)}\kk(\t\times\Aone)$),
and $M_\alpha=
\sum_\lambda^{v\in\ _\lambda V^{f_\alpha}}\O(\t\times\Aone)p_\alpha(v)$
contains $(\pi,\Id,\Id)^*H^\bullet_{\GO\rtimes\Gm}(\Gr_G,S(V))$. 
It follows that the two $W$-actions on
$\left(\bigoplus_\lambda(\Id,\pi,\Id)_*\O(\Gamma_\lambda)\otimes\ _\lambda V
\right)\otimes_{\O(\t\times\Aone)}\kk(\t\times\Aone)$ coincide.
In particular, we have defined unambiguously a $W$-action on
$\left(\bigoplus_\lambda(\Id,\pi,\Id)_*\O(\Gamma_\lambda)\otimes\ _\lambda V
\right)\otimes_{\O(\t\times\Aone)}\kk(\t\times\Aone)$.

Now we claim that $(\pi,\Id,\Id)^*\phi(V)=
\bigcap_{w\in W}^{\alpha\operatorname{simple}}w(M_\alpha)=
(\pi,\Id,\Id)^*H^\bullet_{\GO\rtimes\Gm}(\Gr_G,S(V))$. In effect, 
note that if $(\Id,\pi,\Id)\Gamma_\nu\cap(\Id,\pi,\Id)\Gamma_\mu$
has codimension 1 in $(\Id,\pi,\Id)\Gamma_\nu$ (and in
$(\Id,\pi,\Id)\Gamma_\mu$ as well), then necessarily $\mu=\nu+k\beta$
for certain integer $k$, and certain positive root $\beta$ (not necessarily
a simple root). Let us choose $w\in W$ such that $\alpha:=w^{-1}(\beta)$
is a simple root. Then $w^{-1}(\mu)=w^{-1}(\nu)+k\alpha$.
Since we know that any section in $M_\alpha$ extends through the generic
point of $(\Id,\pi,\Id)\Gamma_{w^{-1}(\mu)}\cap(\Id,\pi,\Id)
\Gamma_{w^{-1}(\mu)+k\alpha}$, we conclude that any section in
$w(M_\alpha)$ extends through the generic point of 
$(\Id,\pi,\Id)\Gamma_\mu\cap(\Id,\pi,\Id)\Gamma_\nu$.
It follows that any section in 
$\bigcap_{w\in W}^{\alpha\operatorname{simple}}w(M_\alpha)$ is regular
off a codimension 2 subvariety. Since both $(\pi,\Id,\Id)^*\phi(V)$ and
$(\pi,\Id,\Id)^*H^\bullet_{\GO\rtimes\Gm}(\Gr_G,S(V))$ are flat
$\O(\t\times\Aone)$-modules coinciding with 
$\bigcap_{w\in W}^{\alpha\operatorname{simple}}w(M_\alpha)$
generically, we conclude that 
$(\pi,\Id,\Id)^*\phi(V)=
\bigcap_{w\in W}^{\alpha\operatorname{simple}}w(M_\alpha)=
(\pi,\Id,\Id)^*H^\bullet_{\GO\rtimes\Gm}(\Gr_G,S(V))$.
This completes the proof of the a).

b) follows by comparing Propositions~\ref{tenstop} and~\ref{tensalg}.
\hfill $\Box$

\subsection{Cohomology is fully faithful}
The following property of the cohomology functor will play an
important role in the proof of the main results.
Let $\CatTopss$ (resp. $\widetilde{\mathcal{IC}}$) 
denote the full subcategory of 
semisimple complexes in $D_{\GO\rtimes \Gm}(\Gr)$
(resp. $D_\GO(\Gr)$).

\begin{lem}
\label{fullyTop}
a) The functor $H_{\GO\rtimes\Gm}^\bullet:\
\CatTopss\to Coh^\Gm(N_{(\Lt^*/W)^2}\Delta)$ is a full
imbedding.

b) The functor $H_\GO^\bullet:\ 
\widetilde{\mathcal{IC}}\to Coh^\Gm(\bT(\Lt^*/W))$
is a full imbedding.
\end{lem}

\proof is due to V.~Ginzburg, see~\cite{G0}.
We prove a), and the proof of b) is identical.
For $V_1,V_2\in Rep(\LG)$ we have
$\Ext_{\GO\rtimes\Gm}^\bullet(S(V_1),S(V_2))=
\Ext_{G\times\Gm}^\bullet(S(V_1),S(V_2))$. Let us denote by $Res_G^T$
the forgetting functor from the $G\times\Gm$-equivariant derived
category to the $T\times\Gm$-equivariant derived category. Then the
Weyl group $W$ acts naturally on 
$\Ext_{T\times\Gm}^\bullet(Res_G^TS(V_1),Res_G^T(V_2))$, and 
$\Ext_{G\times\Gm}^\bullet(S(V_1),S(V_2))=
\Ext_{T\times\Gm}^\bullet(Res_G^TS(V_1),Res_G^T(V_2))^W$.
On the other hand, we know by Theorem~\ref{compaThm}~a) that
$H^\bullet_{T\times\Gm}(\Gr,S(V_{1,2}))\simeq\varphi(V_{1,2})$, and
clearly, $\Hom(\phi(V_1),\phi(V_2))=\Hom(\varphi(V_1),\varphi(V_2))^W$.
Thus it suffices to prove that
$$\Ext_{T\times\Gm}^\bullet(Res_G^TS(V_1),Res_G^T(V_2))\simeq$$
$$\Hom_{\CO(\t\times(\t/W)\times\Aone)}
\left(H^\bullet_{T\times\Gm}(\Gr,S(V_1)),
H^\bullet_{T\times\Gm}(\Gr,S(V_2))\right)=$$
$$\Hom_{H^\bullet_{T\times\Gm}(\Gr)}\left(H^\bullet_{T\times\Gm}(\Gr,S(V_1)),
H^\bullet_{T\times\Gm}(\Gr,S(V_2))\right).$$

Following~\cite{BL}, recall the definition of the LHS: we choose
finite dimensional approximations $P_i$ to the classifying space of 
$T\times\Gm$, and we have ind-varieties $P_i\Gr$ fibered over $P_i$
with fibers isomorphic to $\Gr$. We also have semisimple perverse
sheaves $P_iS(V_{1,2})$ on $P_i\Gr$, and finally
$\Ext_{T\times\Gm}^\bullet(Res_G^TS(V_1),Res_G^T(V_2))$ is defined as
$\underset{\to}{\lim}\Ext_{P_i\Gr}^\bullet(P_iS(V_1),P_iS(V_2))$.
Since $T\times\Gm$ is a torus, we can choose $P_i$ to be the products
of projective spaces (of increasing dimension).
We can choose a generic action of $\Gm$ on $P_i\Gr$ (linear along
$P_i$, and via a one-parametric subgroup of $T\times\Gm$ along $\Gr$)
such that the corresponding Bialynicki-Birula decomposition of
$P_i\Gr$ is cellular. Then we can apply the Theorem of~\cite{G0} to
conclude that $\Ext_{P_i\Gr}^\bullet(P_iS(V_1),P_iS(V_2))\simeq
\Hom_{H^\bullet(P_i\Gr)}(H^\bullet(P_i\Gr,S(V_1)),H^\bullet(P_i\Gr,S(V_2)))$. 
But the limit of the RHS as $i$ grows is 
$\Hom_{H^\bullet_{T\times\Gm}(\Gr)}\left(H^\bullet_{T\times\Gm}(\Gr,S(V_1)),
H^\bullet_{T\times\Gm}(\Gr,S(V_2))\right)$.

The lemma is proved.
\hfill $\Box$

\subsection{Proof of Theorems~\ref{free_cohh},~\ref{free_coh}}
The Theorems follow from Theorem~\ref{compaThm} in view of 
Lemmas~\ref{propKons}(b),~\ref{fullyTop}.

\subsection{Homology}
The goal of this subsection is to express equivariant cohomology of
arbitrary (not necessarily semisimple) equivariant complexes in terms
of Harish-Chandra bimodules, and to prove Theorem~\ref{TodaTh}.
By Theorem~\ref{free_cohh}
we  have  a monoidal equivalence $\eqss$ 
between the category $\HC_\hbar^{fr}$ of free asymptotic
Harish-Chandra bimodules and the category $\CatTopss$ of 
semisimple complexes in $D_{\GO\rtimes \Gm}(\Gr)$.

A standard argument shows any functor on free Harish-Chandra bimodules
is representable by a unique up to a unique isomorphism
(not necessarily free) Harish-Chandra bimodule. Thus we have 
a functor from $\fufu:D_{\GO\rtimes \Gm}(\Gr)\to \widetilde \HCh$
equipped
with a functorial isomorphism
$$\Hom(\eqss(F), \F)\cong \Hom(F,\fufu(\F));$$
$\fufu$ is defined uniquely up to a unique isomorphism.

Since $\eqss$ is a full embedding,
 $\fufu \circ \eqss\cong\Id$ canonically, i.e.
$\fufu$ restricted to the category of semi-simple complexes is the inverse
equivalence to $\eqss$.
It is easy to see from the definition that $\fufu$ is a homological functor;
thus it actually lands in the category $\HCh$ of finitely generated
Harish-Chandra bimodules,
$\fufu: D_{\GO\rtimes \Gm}(\Gr)\to \HCh$.

We will say that  a functor $F$ between two monoidal
categories is quasi-monoidal if a functorial map 
$F(X)*F(Y)\to F(X*Y)$ is fixed for any objects $X,Y$ of the target category,
compatible with the associativity isomorphisms in the two categories in 
the natural sense.

\begin{prop}\label{fufux}
a) $\fufu$ carries a unique quasi-monoidal structure, whose
restriction to $\CatTopss$ induces the natural monoidal
structure on $\Id_{\Bimodfr}\cong \fufu \circ \eqss$.

b) We have a natural isomorphism 
$H^\bullet_{\GO\rtimes \Gm}\cong \qkap \circ \fufu$
compatible with the (quasi)monoidal structure.

\end{prop}

\proof a) 
It is easy to see that 
$$\fufu (\F)= R\Hom(IC_0, \O(\LG)* \F)\cong R\Hom(IC_0,
\F*\O(\LG))$$
canonically
(where the last isomorphism comes from the fact that
both sides can be identified with $\oplusl_\lambda 
RHom(V_\lambda^*\otimes IC_\lambda, \F)$). Then 
the quasi-monoidal structure comes from the coalgebra structure on 
$\O(\LG)$.
Uniqueness follows from the fact that for every $\F\in D_{\GO\rtimes \Gm}$
we can find a free asymptotic Harish-Chandra bimodule $V$ and a surjection
$V\to \fufu (\F)$; by definition of $\fufu$ it comes from a map
 $L=\eqss(V)\to \F$.
Then, in view of functoriality, for $\F_1,\,\F_2\in D_{\GO\rtimes \Gm}$
 the map $\fufu(\F_1)*\fufu(\F_2)\to
\fufu(\F_1*\F_2)$ is uniquely determined by the isomorphism
 $\fufu(L_1)*\fufu(L_2)
\to \fufu(L_1*L_2)$ for $L_i\to \F_i$ ($i=1,2$) as above.

b) For $\F\in D_{\GO\rtimes \Gm}(\Gr)$ we can find an exact sequence
$V_1\to V_2\to \fufu(\F)\to 0$, where $V_i$ are 
free asymptotic Harish-Chandra bimodules. 
To this sequence there corresponds a sequence
of maps in $D_{\GO\rtimes \Gm}(\Gr)$ with zero composition:
 $L_1\to L_2\to \F$, where $L_i=\eqss(V_i)$.
We have $\qkap\circ \fufu (\F) = CoKer (\qkap(V_1) \to \qkap(V_2))=CoKer
(H^\bullet_{\GO\rtimes \Gm}(L_1)\to 
H^\bullet_{\GO\rtimes \Gm}(L_2))$. The latter module
 maps canonically
to $H^\bullet_{\GO\rtimes \Gm}(\F)$. Thus we defined a map $\qkap \circ 
\fufu (\F)\to H^\bullet_{\GO\rtimes \Gm}
(\F)$. A standard argument shows that this map does not depend on the 
choice of $V_1,\, V_2$.
We have obtained a natural transformation between the two functors. This
 transformation is an isomorphism
on semi-simple complexes. Since both functors are homological (where
we use exactness of $\qkap$, Lemma \ref{propKons}(a)) and semi-simple objects generate the
 triangulated category
$D_{\GO\rtimes \Gm}(\Gr)$, we see that the transformation 
is an isomorphism. \hfill $\Box$

We can clearly extend all of the above to  Ind-objects.
We will be particularly interested in the Ind-object
$\ds=\varinjlim \ds_\lambda$, where $\ds_\lambda$ is the dualizing
sheaf of the closure of the $\GO$-orbit $\Gr_\lambda$.
Notice that $H_\bu^{\GO\rtimes \Gm}(\Gr)=H^\bu_{\GO\rtimes \Gm}(\ds)$,
 and the convolution algebra
structure on  homology comes from the structure of an algebra in the 
monoidal category $D_{\GO\rtimes \Gm}$ on the object $\ds$ and the monoidal
structure on the cohomology functor.% {\bf does this need a proof?}.

We will now describe the corresponding object in the category of Harish-Chandra
bimodules. Recall first the duality $\prd$ on $\HChfr$,
$\prd:M\mapsto \Hom_{U_\hbar}(M,U_\hbar)$, where $\Hom$ is taken
with respect to the right action of
$U_\hbar$, and the left action on $M,\, U_\hbar$ is used to define,
respectively, the right and the left action on the $\Hom$ module.
Let $Forg$ denote the forgetful functor from Harish-Chandra bimodules
to vector spaces. The functor $Forg\circ \prd$ is represented
by the Ind-object, which is readily identified with $U_\hbar \otimes
\O(\LG)$, the module of $\hbar$-differential operators on $\LG$.
Furthermore,  $U_\hbar \otimes
\O(\LG)$ carries another commuting structure of a Harish-Chandra bimodule,
and we set $\cK=\psf(U_\hbar\otimes \O(\LG))$, where this second structure
of a  Harish-Chandra bimodule is used to compute $\psf$. 
 Thus
$\cK$ is an Ind-object in the category of Harish-Chandra bimodules;
moreover, it is an algebra Ind-object in this monoidal category.

\begin{prop}\label{fufudsK}
 We have a canonical isomorphism of algebra Ind-objects:
$\fufu(\ds)\cong \cK$.
\end{prop}

To prove Proposition we need another auxiliary

\begin{lem}\label{VerD} We have a canonical isomorphism
 $\Ver\circ \eqss\cong \eqss\circ   \inv_{\LG}\circ \prd$, where
$\Ver$, %$\prd$,
 $\inv_{\LG}$ denote, respectively, Verdier duality
%, projective module duality $M\mapsto Hom_{U_\hbar}(M,U_\hbar)$ 
and 
the functor induced by the canonical outer automorphism of $\LG$
interchanging conjugacy classes of $g$ and $g^{-1}$, $g\in \LG$
(the Chevalley involution).
\end{lem}

\proof: Recall that a monoidal category ${\mathcal C}$ is rigid 
if for any object $V\in {\mathcal C}$
there exists another object $V^*\in {\mathcal C}$ 
and morphisms $\iota:{\bf 1}\to V\otimes V^*$ and $\tau:V^*\otimes V\to 
{\bf 1}$ satisfying a certain compatibility constraint, see~\cite{DM}. 
Given $V$, an object $V^*$ together with morphisms $\iota,\tau$ is unique
up to a unique isomorphism if it exists. Thus for a rigid
category  ${{\mathcal C}}$ we have a canonical (up to a canonical isomorphism)
functor $ {{\mathcal C}}\mapsto  {{\mathcal C}}^{op}$ sending $V$ to $V^*$.
It is immediate to check that monoidal categories $\PerSatake$
and $\HC_\hbar^{fr}$ are rigid categories with duality functors given by
$\inv_G\circ \Ver$ and $\prd$ respectively,
where the functor $\inv_G$ is induced by the Chevalley involution 
of $G$. The equivalence between the two
categories intertwines the canonically defined dualities.
Also, it is well-known that $\eqss \circ \inv_{\LG}\cong \inv _G\circ
\eqss$ canonically. The Lemma follows. 
\hfill $\Box$

{\em Proof of  Proposition \ref{fufudsK}.} The Ind-object
$\ds$ represents the functor $\F\mapsto H^\bullet_{\GO\rtimes \Gm}(\Ver(\F))$.
In view of Lemma \ref{VerD} and Proposition \ref{fufux}(c)
we see that 
 the Ind-object $\fufu(\ds)$ represents the functor
$M\mapsto \psf(\prd (M))$ on the category of free asymptotic
Harish-Chandra modules. 
 It is straighforward to see from the definitions
that the Ind-object $\cK$ represents the same functor.
The isomorphism of functors yields an isomorphism of Harish-Chandra
bimodules. Since the isomorphism of functors is compatible with the
monoidal structure, the isomorphism of Harish-Chandra bimodules
is compatible with the algebra structure.
%The claim follows {\bf though it may be nice to say something about
%the algebra structure}. 
\hfill $\Box$

{\em Proof of Theorem \ref{TodaTh}.}
By Propositions \ref{fufux} and \ref{fufudsK} we have
an isomorphism of algebras
$$H^{\GO\rtimes \Gm}_\bu(\Gr)=H_{\GO\rtimes \Gm}^\bu(\ds)\cong
\psf(\cK).$$
The latter is by definition the algebra of the quantum Toda
lattice. \hfill $\Box$

\subsection{Formality from purity}
In this section we combine the above Ext computation with a standard argument
which allows one to derive formality of the RHom algebra from purity of the Ext
spaces. The result is a description of the derived Satake category, including
the version equivariant with respect to the loop rotation. 

{\em Except for some technical details, this section does not contain
original contributions of the authors.} We have learned the
geometric (respectively, algebraic) ideas exposed
here from V.~Ginzburg (respectively, L.~Positselski) around 1998.

In order to be able to use Frobenius weights we extend the basic setting,
and consider $\GO$, $\Gr$ etc. over $\Fqbar$, and the categories of
equivariant $l$-adic sheaves on $\Gr$. 
%All the above arguments apply in this setting mutatis mutandis. 

Consider the following general situation. Let $R$ be a finitely localized
ring of integers of a number field $E$, 
and let $\Fq$ be a finite field quotient of $R$. Let $X_R$ be a flat scheme
over $R$ acted upon by a smooth affine group scheme $G_R$, such that the set of
orbits is finite. 
%More generally, for application to affine Grassmannian,
%we may take $X_R$ to be an ind-scheme,
%and $G_R$ a proalgebraic group scheme, such that in all the terms of the 
%inductive sistem the sets of orbits are finite. 
We denote by 
$(X_{\Fqbar},G_{\Fqbar})$ (resp. $(X_{\bar E},G_{\bar E})$) the base change
of $(X_R,G_R)$ to a geometric point of $R$ over $\Fq$ (resp. over the 
generic point). We choose a prime $l$ invertible in $R$.
Let $D_{G_{\Fqbar}}(X_{\Fqbar})$ (resp. $D_{G_{\bar E}}(X_{\bar E})$) stand 
for the bounded equivariant constructible derived category of 
{\'e}tale $\Qlbar$-sheaves on $X_{\Fqbar}$ (resp. $X_{\bar E}$)
(see e.g.~\cite{BD},~7.4).
We choose an isomorphism $\Qlbar\simeq\kk$
(under a technical assumption that $\kk$ has the same cardinality 
as $\Qlbar$), and an embedding 
${\bar E}\hookrightarrow\BC$, and we denote by $(X_\BC,G_\BC)$ the base
change of $(X_{\bar E},G_{\bar E})$ to $\BC$. 
Let $D_{G_\BC}(X_\BC)$ (resp. $D_{G_\BC}^{top}(X_\BC,\Qlbar),\
D_{G_\BC}^{top}(X_\BC)$) stand for the bounded equivariant constructible
derived category of {\'e}tale $\Qlbar$-sheaves on $X_\BC$ (resp.
bounded equivariant constructible derived category of sheaves with
$\Qlbar$-coefficients, resp. $\kk$-coefficients, 
{\em in the classical topology} of $X_\BC$).

\begin{prop}
\label{bbdd}
There exists a localization $R_{(r)}$ of $R$ such that for any point
$R_{(r)}\twoheadrightarrow\Fq$, we have the following 
chain of natural equivalences:
$$D_{G_{\Fqbar}}(X_{\Fqbar})\stackrel{\alpha}{\to}
D_{G_{\bar E}}(X_{\bar E})\stackrel{\beta}{\to}
D_{G_\BC}(X_\BC)\stackrel{\gamma}{\to}
D_{G_\BC}^{top}(X_\BC,\Qlbar)\stackrel{\delta}{\to}
D_{G_\BC}^{top}(X_\BC)$$
\end{prop}

{\em Sketch of proof.} The argument is taken from~\cite{BBD},~6.1.
The first equivalence $\alpha$ is constructed
in~\cite{BBD},~6.1.9 (existence of good
models). To justify the finiteness assumptions of {\em loc. cit.} we note
that the set of isomorphism classes of $G$-equivariant irreducible perverse
sheaves on $X$ is finite. Since the equivariant derived categories are not
considered in {\em loc. cit.} we note that according to~\cite{BD},~7.4,
to compare the equivariant Exts between equivariant irreducible perverse 
sheaves on $X$ it suffices to compare the usual Exts between the lifts of
these sheaves to $X\times G^n$ (see the canonical spectral sequence~(312) of
{\em loc. cit.}). These are calculated by the K{\"u}nneth formula.

The second equivalence $\beta$ is just the base change from $\bar E$ to
$\BC$. The third equivalence $\gamma$ is the classical M.~Artin's comparison
theorem of {\'e}tale and classical cohomology, see~\cite{BBD},~6.1.2(B$''$). 
Finally, $\delta$ is induced by our isomorphism $\Qlbar\simeq\kk$.
\hfill $\Box$

\bigskip

From now on we will restrict our attention to the equivariant derived
categories in our new setting, that is, over $\Fqbar$.

Let $X$ be an algebraic variety over a finite $\Fq$, and
let $D(X_{\Fqbar})$ stand for the bounded constructible $l$-adic
derived category of $X_{\Fqbar}$.
Let $\F$ be a pure weight zero object of the 
$l$-adic derived category of $X$ of geometric origin.
 The space $\Ext^i(\F_{\Fqbar},\F_{\Fqbar})=\oplusl
\Ext^i_j(\F_{\Fqbar},\F_{\Fqbar})$ carries a canonical
grading by Frobenius weights;
here the subindex
denotes base change to $\Fqbar$, and $\Ext^i_j$ is the component 
of weight $j$. Recall that by Deligne Theorem~\cite{W2}
 $\Ext^i_j=0$ for $j<i$.

Let $\Epur$ be a graded algebra and $\phi:\Epur\to \Ext^\bu(\F_{\Fqbar},\F_{\Fqbar})^{op}$ 
be a homomorphism
sending a graded component $\Epur^i$ to $Ext^i_i(\F_{\Fqbar},\F_{\Fqbar})$.
%$\Epur=\oplusl Ext^i_i(\F_{\Fqbar},\F_{\Fqbar})$, and let us consider $\Epur$ as a dg-algebra
%with zero differential. 

We will consider the graded algebra $\Epur$ as a dg-algebra with zero differential.

\begin{prop}\label{general}
There exists a canonical functor 
$\Phi_X:\ D_{perf}(\Epur)\to D(X_{\Fqbar})$ sending the free module
to $\F_\Fqbar$ and inducing the map $\phi$ on Ext groups. 
\end{prop}

{\em Sketch of proof.}  The complex $\F_\Fqbar$ is semi-simple, i.e.\ is isomorphic
to $\oplus \IC_i[d_i]$, where $\IC_i$ is an irreducible perverse sheaf and $d_i\in \Zet$.
By Beilinson's Theorem \cite{Bedr} the $l$-adic derived 
category
contains the derived category of perverse sheaves as a full subcategory, thus 
$\Ext^\bu(\F_{\Fqbar},\F_{\Fqbar})$ coincides with Ext in the category of perverse sheaves.
We can assume without loss of generality that $\Epur$ is finitely presented; thus the map $\phi$
factors through a map $\phi_{fin}\colon \Epur\to \Ext^\bu_{\A} (\F_{\Fqbar},\F_{\Fqbar})$, where $\A$ 
is the Serre subcategory in the category of perverse sheaves on $X_\Fqbar$
generated by a finite set of irreducible 
objects, including $\IC_i$. Moreover, the argument 
of~\cite{Bedr} (cf. Examples in {\em loc.cit.}, 1.2, p.28) shows
that all irreducible objects of $\A$ can be assumed to be of geometric origin.

We can identify the abelian category $\A$ with the category of finite length $A$-modules,
where the pro-finite dimensional algebra (algebra in the 
tensor category of pro-finite
dimensional vector spaces) $A$ is defined by
$A=\End(\oplus \P_s)$; here $s$ runs over the (finite) set of isomorphism
classes of irreducible objects in $\A$, and  $\P_s$ is a pro-object in 
$\A$ which is a projective
cover of the corresponding irreducible object $L_s$, cf.~\cite{B}.
We fix an isomorphism $Fr_q^*(L_s)\cong L_s$, which induces a pure weight 
zero Weil structure
on $L_s$ (this is possible because $L_s$ 
has geometric origin, see~\cite{BBD}).
Since projective cover of an irreducible object is unique 
up to a non-unique isomorphism,
we can (and will) fix an isomorphism $Fr_q^*(\P_s)\cong \P_s$. 
Then conjugation with Frobenius
is an automorphism of $A$ (which we will also call Frobenius).

By a result of~\cite{BBD}, Frobenius acts on $\Ext^1(L_s,L_{s'})$ 
with positive weights. It follows
that Frobenius finite elements are dense in $A$, and they form a graded subalgebra $A^{gr}$ with
finite dimensional graded components, where
the grading comes from Frobenius weights. Moreover, components of negative degree in $A^{gr}$ vanish,
while $A^{gr}_0$ is semisimple.
Obviously, $\A$ is identified with the category of finite length $A^{gr}$ modules, on which $A^{gr}_N$
acts by zero for $N\gg 0$.

We now consider the object 
 $L=\oplus L_i[d_i]\in D^b(\A)$ (where $L_i$ corresponds to $\IC_i$) and a dg-algebra $\sD=\RHom_\A(L,L)$
(well defined as an object of the category of dg-algebras with inverted quasi-isomorphisms).
Recall that we have an equivalence $D^b(\A)\cong D_{perf}(\sD^{op})$, $M\mapsto \RHom(L,M)$.
We lift $L$ to an object $\tilde L=\oplus \tilde L_i[d_i](d_i)$ of the derived category of graded 
$A^{gr}$-modules, where $\tilde L_i$ is the irreducible $A$-module concentrated in degree zero,
and $(d)$ stands for shift of grading by $d$. Then the algebra $\Ext^\bu(L,L)$ acquires
an additional grading, and $\sD$ can be chosen to carry also an additional grading compatible
with the grading on $\Ext$'s. 
We have a homomorphism $\Epur\to\oplusl_i\Ext^i_i(L,L)$, where
the lower index denotes the additional grading, 
and the fact that $A^{gr}$ is positively graded
implies that $\Ext^i_j(L,L)=0$ for $j<0$. 
Thus existence of a canonically defined functor $D_{perf}(\Epur)\to D^b(\A)$ follows from the 
standard Lemma \ref{gradualg}a) below. 

We leave it as an exercise to the reader to show that the composed
functor $\Phi_X:\ D_{perf}(\Epur)
\to D^b(\A)\to D(X_{\Fqbar})$ does not depend on the choice of $\A$ 
up to a canonical isomorphism. \hfill $\Box$

We will also have to use functoriality properties of the above 
construction. We spell these out now.

\begin{prop}\label{foon}
a) Let $X$, $Y$ be algebraic varieties over a finite field $\Fq$, and 
%$F:D(X) \to D(Y)$,
$F:D(X_{\Fqbar})
\to D(Y_{\Fqbar})$ be a functor
 satisfying the following conditions

\begin{enumerate} \item $F$ commutes with the pull back under Frobenius
functor, i.e. an isomorphism $F\circ Fr_X\cong Fr_Y \circ F$ is fixed.
%The pair $(F,\tilde F)$ is compatible with
%base change functors $D(X_{\Fqbar})\to D(X)$, $D(Y_{\Fqbar})\to D(Y)$.

\item \label{purprop}
 $F$ sends pure weight zero Weil complexes to pure weight zero Weil complexes.

\item \label{exaprop}
Let $\A$ be a finitely generated Serre subcategory in 
$Perv(X_{\Fqbar})$ invariant under the Frobenius pull-back functor.
%, and $\tilde \A$ be its preimage
%under the base change functor  $D(X_{\Fqbar})\to D(X)$.
 Then there exists a
natural exact functor $F_{\A}:Com(\A)\to Com(Perv(Y_{\Fqbar}))$,
%, $F_{\tilde \A}:Com(\tilde \A) \to Com(Perv(Y_{\Fqbar}))$, 
compatible with Frobenius and equipped with a natural isomorphism
$\beta_Y\circ F_{\A}\iso F\circ\beta_X$. Here $\beta_X:\ Com(\A)\to
D(X_{\Fqbar}),\ \beta_Y:\ Com(Perv(Y_{\Fqbar}))\to D(Y_{\Fqbar})$ are
the natural functors.
\end{enumerate}

Then for $\F$, $\Epur$ as in Proposition \ref{general}, the construction of Propositon \ref{general}
 is compatible with
$F$, that is, there is a natural isomorphism $\psi_{X\to Y}:\ 
\Phi_Y\iso F\circ\Phi_X$.

b) Assume furthermore that  $ F_1:\ D(Y_{\Fqbar})\to D(Z_{\Fqbar})$  
 is a functor satisfying the above conditions. Then 
 $F_1\circ F$ also satisfies these conditions,
and the two isomorphisms $\psi_{X\to Z},\ F_1\circ\psi_{X\to Y}$ 
between the two functors 
$\Phi_Z\iso F_1\circ F\circ\Phi_X:\  
D_{perf}(\Epur)\to D(Z_{\Fqbar})$ coincide.

\end{prop}

{\em Proof.\ } Let $\CalC\in D_{perf}(\Epur)$.
As above, we find an abelian subcategory $\A'\subset Perv(Y)$ containing all subquotients of the terms of complexes
$F_{\A}(\G)$, $\G\in \A$;  a complex $\CalC'$ of pro-objects in $Perv_{mix}(Y)$ quasiisomorphic
to $F(\CalC)$, whose terms with forgotten Frobenius action are projective pro-objects
in $\A'$. We  also have a dg-algebra $\sD'$ equipped with an additional grading, which acts
on $\CalC'$ in a way compatible with the grading by Frobenius weights, so that the action
induces a quasiisomorphism $\sD'\to \RHom_{\A'}(F(\CalC),F(\CalC))$.

Consider a dg-module $B_{X,Y}$ of $\sD'\otimes \sD^{op}$-module defined by $B_{X,Y}:=
\Hom^\bu(\CalC',F(\CalC))$. It is not hard to see that the composed functor
$D_{perf}(\sD)\cong D^b(\A)\overset{F}{\To} D^b(\A')\cong D_{perf}((\sD')^{op})$ arises from this
bimodule as described in Lemma \ref{gradualg}(b). The action of Frobenius endows $B_{X,Y}$
with an additional grading compatible with the gradings on $\sD,\sD'$. 
Thus  Lemma \ref{gradualg}(b)
provides the sought for isomorphism between the two functors $D_{perf}(\Epur)\to D^b(\A')$.

Finally part (b) of the Proposition 
%
%, given another smooth map $f':Z\to Y$ 
%or a proper map  $g':Y\to Z$ one can 
%deduce an equality of the two isomorphisms between the two functors $D_{perf}(\Epur)
%\to D^b(A'')$, $\A''\subset Perv(Z)$ 
can be deduced from Lemma \ref{gradualg}(c). \hfill $\Box$

\begin{Rem}\label{pull_push_fun}
We will apply Proposition \ref{foon} when
$F=f^*$ or $F=g_*$ for a 
smooth map $f:Y\to X$ or a proper map $g:X\to Y$. Each of  the functors
$F=f^*$, $F=g_*$ 
satisfies the requirements of the Propositions:
 for $f^*$ this is standard, and for $g_*$ 
property \ref{exaprop}  
follows from a construction described in \cite{Bedr}, page~41, 
and \ref{purprop} follows from \cite{BBD}.

Thus Proposition \ref{foon} implies functoriality 
of the construction of Proposition
\ref{general} with respect to proper push-forward and 
smooth pull backs. Also, Proposition
\ref{foon}(b) implies 
compatibility of the isomorphism of Proposition \ref{general} with
the base change isomorphism for 
 a proper map $X\to Y$ and a smooth map $Y'\to Y$. 
\end{Rem}

\begin{Rem}\label{multifoon}
A result similar to Proposition \ref{foon} holds, with a similar proof, for
 a functor
 $F:D(\ ^1X_{\Fqbar})\times \cdots \times D(\ ^nX_{\Fqbar}) \to 
D(Y_{\Fqbar})$. Examples of this situation arise when
$Y=\ ^1X\times \cdots \times\ ^nX$, 
$F:(\F_1, \dots, \F_n) \mapsto \F_1\boxtimes \cdots \boxtimes
\F_n$; or in a twisted version of this situation (see below). 
\end{Rem}

\begin{lem}\label{gradualg}
a) Let $\sD=\oplusl \sD^i_j$ be a dg-algebra equipped with an additional ``inner'' grading
(denoted by a subindex), which is compatible with the differential (thus we have $d:\sD^i_j\mapsto
\sD^{i+1}_j$). Assume that $H^i_j(\sD)=0$ for $j<i$. Then there is a canonical morphism
in the category of dg-algebras with inverted quasiisomorphisms: $H_{pur}=\oplusl_i H^i_i(\sD) \to \sD$.
In particular, we have a canonical push-forward functor $D_{perf}(H_{pur})\to D_{perf}(\sD)$.

b) Let $\sD$, $\sD'$ be dg-algebras satisfying the assumptions 
of (a). Let $B\in D(\sD'\otimes \sD^{op})$ be such that the forgetful functor $D(\sD'\otimes \sD^{op})
\to D(\sD')$ sends $B$ to the free rank one module over $\sD'$; thus $B$ defines
a homomorphism $\phi_B:H^\bu(\sD)\to H^\bu(\sD')$. 
Let $\Epur$, $\phi:\Epur\to \oplusl_iH_i^i(\sD)$
be as above, and consider the functors $\Phi:D_{perf}(\Epur)\to D_{perf}(\sD)$,
$\Phi':D_{perf}(\Epur)\to D_{perf}(\sD')$ arising from $\phi$, $\phi_B\circ \phi$ respectively
by the construction of part(a).

 Consider the functor $\Phi_B:D_{perf}(\sD)\to D_{perf}(\sD')$ given by $M\mapsto B\Lotimes_\sD M$.
Assume that $B$ carries an additional grading compatible with the gradings on $\sD$, $\sD'$.
Then we have a natural isomorphism $\Phi'\cong \Phi_B\circ \Phi$.

c) Let $\sD$, $\sD'$, $\sD''$ be three dg-algebras as above, and $B$, $B''$ be modules
for $\sD'\otimes \sD^{op}$, $\sD''\otimes (\sD')^{op}$ as above. For a homomorphism
$\Epur \to \oplusl H^i_i(\sD)$ the two isomorphisms between the two functors 
$D_{perf}(\Epur)\to D_{perf}(\sD'')$ arising from part (b) coincide.
\end{lem}

{\em Proof.\ } We remind the idea of the construction in part (a), and leave (b,c) to the interested
reader. Let $\sD_i\subset \sD$ be the subcomplex of elements of inner
degree $i$. We have a sub dg-algebra 
$\sD_{up}:= \oplusl_i \tau_{\leq i} \sD_i$, where
we use the standard notation $\tau$ for truncation of a complex.
Furthermore, $\sD_{up}$ has a quotient algebra with zero differential
$\sD_{diag}:=\oplusl_i \tau_{\geq i}\tau_{\leq i} \sD_i$. 
The conditions of part (a) guarantee
that the projection homomorphism $\sD_{up}\to \sD_{diag}$ is a quasi-isomorphism. The composition
of the formal inverse to this quasi-isomorphism 
with the embedding $\sD_{up}\imbed \sD$
is the desired morphism. \hfill $\Box$

\subsection{Proof of  Theorem~\ref{polkovnik}} 
We construct the first equivalence, the second one is similar. 
We will construct a monoidal functor 
$\Psi:D^{\LG}_{perf}(\Duh) \to D_{\GO\ltimes \Gm}(\Gr)$, whose restriction to
the full subcategory $\HChfr\subset D^{\LG}_{perf}(\Duh)$
is identified with $\eqss$ (where the full embedding 
sends a $\LG$-equivariant graded $U_\hbar$-module to the same module
considered as a dg-module with zero differential). Then $\Psi$
sends a set of generators of the source triangulated category
to generators of the target category, and induces an isomorphism on
Hom's between the generators, hence it is an equivalence.

It suffices to construct a collection of functors
$\Psi_\la:  D^{\LG}_{perf}(\Duh)_{\leq \la}\to
D_{\GO\ltimes \Gm}(\Gr_{\leq \la})$, where $\la$ is a coweight of $G$,
$\Gr_{\leq \la}$ is the closure of the corresponding $\GO$ orbit on
$\Gr$, and $ D^{\LG}_{perf}(\Duh)_{\leq \la}$
 is the full subcategory
in  $ D^{\LG}_{perf}(\Duh)$ generated by the objects $V\otimes \Duh$,
where $V$ is an irreducible representation of $\LG$ with a highest weight
$\mu \leq \la$.
These functors will be compatible for comparable coweights
(i.e. we have isomorphisms $\Psi_\mu\cong
 \Psi_\la|_{ D^{\LG}_{perf}(\Duh)_{\leq \mu}}$ for $\mu\leq \la$,
satisfying the obvious compatibility for a triple of coweights
$\nu\leq \mu \leq \la$).

The action of $\GO\ltimes \Gm$ on $\Gr_{\leq \la}$ factors
through a finite dimensional algebraic group $H_\la$, and 
 $D_{\GO\ltimes \Gm}(\Gr_{\leq \la})\cong 
D_H(\Gr_{\leq \la})$ naturally.
To describe $\Psi_\la$ we need to provide the following data:
for a smooth $H_\la$-equivariant
 map $X\to \Gr_{\leq \la}$, where $H_\la$ acts on $X$ freely 
we need to provide a functor 
$D^\LG_{perf}(\Duh)_{\leq \la}\to D(X/H_\la)$, compatible with pull-backs
(i.e. a Cartesian section of the category of resolutions of
$\Gr_{\leq \la}/H_\la$ , cf.~\cite{BL},~2.4.3). 
In view of Theorem \ref{free_cohh} we have
a map $\End^\bu (V_{\leq \la}\otimes \Duh) \iso 
\End^\bu (I_{\leq \la})\to End(I_{\leq \la}^X)$;
here $V_{\leq \la}$ is the sum of all irreducible $\LG$-modules
with a highest weight less or equal than $\la$, $I_{\leq \la}
\in D_{\GO\ltimes \Gm}(\Gr_{\leq \la})$ is the sum of
IC sheaves of all $\GO$ orbits in $\Gr_{\leq \la}$, and 
$I_{\leq \la}^X$ is the pull-back of $I_{\leq \la}$ to $X$.
Thus the required functor is given by 
 Proposition \ref{general} in view of 
 purity of equivariant $\Ext$'s
between IC sheaves on the affine Grassmannian,
 see e.g.~\cite{Ga}.
These functors do indeed form a Cartesian
section  in view of
 Proposition~\ref{foon}, cf. Remark~\ref{pull_push_fun}. 
Compatibility between $\Psi_\la$ and $\Psi_\mu$ for $\mu\leq \la$
is left as an exercise for the reader.

In view of Proposition \ref{foon} (cf. Remark \ref{pull_push_fun}), 
a monoidal structure 
for the constructed functor would follow if we show that the functors
from $[D^\LG_{perf}(\Duh)]^2$, $[D^\LG_{perf}(\Duh)]^3$ to the derived category
of sheaves on the convolution space, respectively, 
triple convolution space, given, respectively,
by $(M_1,M_2) \mapsto \Psi(M_1)\boxtimes^{\GO\ltimes \Gm} \Psi(M_2)$,
and by $(M_1,M_2,M_3) \mapsto \Psi(M_1)\boxtimes^{\GO\ltimes \Gm} \Psi(M_2)
\boxtimes^{\GO\ltimes \Gm} \Psi(M_3)$
 (where $\boxtimes^{\GO\ltimes \Gm}$
denotes the twisted external product on the convolution space), are
 compatible with the functors stemming from the construction
of Proposition \ref{general}. This follows from Remark \ref{multifoon}. \hfill $\Box$

\footnotesize{
{\bf R.B.}: Department of Mathematics,
Massachusetts Institute of Technology,
Cambridge, MA 02139, USA;\\ 
%\hphantom{x}\ab\, 
{\tt bezrukav@math.mit.edu}}

\footnotesize{
{\bf M.F.}: Institute of Information Transmission Problems, and 
Independent Moscow University,
Bolshoj Vlasjevskij pereulok, dom 11,
Moscow 119002 Russia;\\
%\hphantom{x}\ab\, 
{\tt fnklberg@gmail.com}}

\end{document}